\documentclass[preprint,11pt]{article}

\usepackage{jmlr2e}
\usepackage{amsmath,amssymb}
\usepackage{array}
\usepackage{cellspace}
\usepackage{multirow}
\usepackage{tabularx}
\usepackage{hyperref}
\usepackage{booktabs}
\newtheorem{Theorem}{Theorem} 
\newtheorem{Lemma}{Lemma} 

\newtheorem{asm}{Assumption}
\newtheorem{Remark}{Remark}
\usepackage{algorithm}

\begin{document}

\title{Entry-Wise Eigenvector Analysis and Improved Rates for Topic Modeling on Short Documents}
\author{\name Zheng Tracy Ke \email zke@fas.harvard.edu \\
       \addr Department of Statistics\\
       Harvard University\\
       Cambridge, MA 02138, USA
       \AND
       \name Jingming Wang \email jingmingwang@fas.harvard.edu \\
       \addr Department of Statistics\\
       Harvard University\\
       Cambridge, MA 02138, USA 
}

\maketitle

\begin{abstract}
Topic modeling is a widely utilized tool in text analysis. We investigate the optimal rate for estimating a topic model. Specifically, we consider a scenario with $n$ documents, a vocabulary of size $p$, and document lengths at the order $N$. When $N\geq c\cdot p$, referred to as the long-document case, the optimal  rate is established in the literature at $\sqrt{p/(Nn)}$. However, when $N=o(p)$, referred to as the short-document case, the optimal rate remains unknown. In this paper, we first provide new entry-wise large-deviation bounds for the empirical singular vectors of a topic model. We then apply these bounds to improve the error rate of a spectral algorithm, Topic-SCORE. Finally, by comparing the improved error rate with the minimax lower bound, we conclude that the optimal rate is still $\sqrt{p/(Nn)}$ in the short-document case.
\end{abstract}


\tableofcontents

\newpage

\section{Introduction}

In today's world, an immense volume of text data is generated in scientific research and in our daily lives. This includes research publications, news articles, posts on social media, electronic health records, and many more. Among the various statistical text models, the topic model \citep{hofmann1999,blei2003latent} stands out as one of the most widely used. 
Given a corpus consisting of $n$ documents written on a vocabulary of $p$ words, let $X=[X_1, X_2,\ldots, X_n]\in\mathbb{R}^{p\times n}$ be the word-document-count matrix, where $X_i(j)$ is the count of the $j$th word in the $i$th document, for $1\leq i\leq n$ and $1\leq j\leq p$. Let $A_1, A_2, \ldots, A_K\in\mathbb{R}^{p}$ be probability mass functions (PMFs). We call each $A_k$ a topic vector, which represents a particular distribution over words in the vocabulary.  For each $1\leq i\leq n$, let $N_i$ denote the length of the $i$th document, and let $w_i\in\mathbb{R}^K$ be a weight vector, where $w_i(k)$ is the fractional weight this document puts on the $k$th topic, for $1\leq k\leq K$. 
In a topic model, the columns of $X$ are independently generated, where the $i$th column satisfies:
\begin{equation} \label{def:model}
X_i \sim  \mathrm{Multinomial}(N_i, \; d_i^0), \qquad \mbox{with}\quad d_i^0 = \sum_{k=1}^K w_i(k)A_k. 
\end{equation}
Here $d^0_i\in\mathbb{R}^p$ is the population word frequency vector for the $i$th document, which admits a convex combination of the $K$ topic vectors. The $N_i$ words in this document are sampled with replacement from the vocabulary using probabilities in $d_i^0$; as a result, the word counts follow a multinomial distribution. Under this model, $\mathbb{E}[X]$ is a rank-$K$ matrix. The statistical problem of interest is using $X$ to estimate the two parameter matrices $A=[A_1, A_2,\ldots, A_K]$ and $W=[w_1, w_2, \ldots, w_n]$.  

Since the topic model implies a low-rank structure behind the data matrix, spectral algorithms \citep{ke2023recent} have been developed for topic model estimation. Topic-SCORE \citep{ke2022using} is the first spectral algorithm in the literature. It conducts singular value decomposition (SVD) on a properly normalized version of $X$, then uses the first $K$ left singular vectors to estimate $A$, and finally uses $\hat{A}$ to estimate $W$ by weighted least-squares. 
 \cite{ke2022using} showed that the error rate on $A$ is $\sqrt{p/(n N)}$ up to a logarithmic factor, where $N$ is the order of document lengths. It matches with the minimax lower bound \citep{ke2022using} when $N\geq c\cdot p$ for a constant $c>0$, referred to as the long-document case. However, there are many application scenarios with $N=o(p)$, referred to as the short-document case. For example, if we consider a corpus consisting of abstracts of academic publications (e.g., see \cite{ke2023recent}), $N$ is usually between 100 and 200, but $p$ can be a few thousands or even larger. 
In this short-document case,  \cite{ke2022using} observed a gap between the minimax lower bound and the error rate of Topic-SCORE. They posted the following questions: Is the optimal rate still $\sqrt{p/(Nn)}$ in the short-document case? If so, can spectral algorithms still achieve this rate? 

In this paper, we give answers to these questions. We discovered that the gap between the minimax lower bound and the error rate of Topic-SCORE in the short-document case came from the unsatisfactory entry-wise large-deviation bounds for the empirical singular vectors. While the analysis in \cite{ke2022using} is effective for long documents, there is considerable room for improvement in the short-document case. We use new analysis to obtain much better large-deviation bounds when $N=o(p)$. Our strategy includes two main components: one is an improved non-stochastic perturbation bound for SVD allowing severe heterogeneity in the population singular vectors, and the other is leveraging a decoupling inequality \citep{de1995decoupling} to control the spectral norm of a random matrix with centered multinomial-distributed columns. 
These new ideas allow us to obtain satisfactory entry-wise large-deviation bounds for empirical singular vectors across the entire regime of $N\geq \log^3(n)$. As a consequence, we are able to significantly improve the error rate of Topic-SCORE in the short-document case. This answers the two questions posted by \cite{ke2022using}:  
The optimal rate is still $\sqrt{p/(Nn)}$ in the short-document case, and Topic-SCORE still achieves this optimal rate. 

Additionally, inspired by our analysis, we have made a modification to Topic-SCORE to better incorporate document lengths.  
We also extend the asymptotic setting in \cite{ke2022using} to a weak-signal regime allowing the $K$ topic vectors to be extremely similar to each other.  

\subsection{Related Literature}

Many topic modeling algorithms have been proposed in the literature, such as LDA \citep{blei2003latent}, the separable NMF approach \cite{Ge,arora2013practical}, the method in \cite{bansal2014provable} that uses a low-rank approximation to the original data matrix, Topic-SCORE \citep{ke2022using}, and LOVE \citep{bing2020fast}. Theoretical guarantees were derived for these methods, but unfortunately, most of them had non-optimal rates even when $N\geq c\cdot p$.
Topic-SCORE and LOVE are the two that achieve the optimal rate when $N\geq c\cdot p$. However, LOVE has no theoretical guarantee when $N=o(p)$; Topic-SCORE has a theoretical guarantee across the entire regime, but the rate obtained by  \cite{ke2022using} is non-optimal when $N=o(p)$. Therefore, our results address a critical gap in the existing literature by determining the optimal rate for the short-document case for the first time. 

Entry-wise eigenvector analysis \citep{erdHos2013spectral,fan2018eigenvector,fan2019simple,abbe2020entrywise,chen2021spectral,ke2022optimal} provides large-deviation bounds or higher-order expansions for individual entries of the leading eigenvectors of a random matrix.  
There are two types of random matrices, i.e., the Wigner type (e.g., in network data and pairwise comparison data) and the Wishart type (e.g., in factor models and spiked covariance models \citep{paul2007}). 
The random matrices in topic models are the Wishart type, and
hence, techniques for the Wigner type, such as the leave-one-out approach \citep{ke2022optimal}, are not a good fit. 
We cannot easily extend the techniques \citep{fan2018eigenvector,chen2021spectral}  for spiked covariance models  either. One reason is that the multinomial distribution has heavier-than-Gaussian tails (especially for short documents), and 
using the existing techniques only give non-sharp bounds. Another reason is the severe word frequency heterogeneity \citep{zipf2013psycho} in natural languages, which calls for bounds whose orders are different for different entries of an eigenvector. Our analysis overcomes these challenges.

\subsection{Organization and Notations}

The rest of this paper is organized as follows. Section~\ref{sec:main} presents our main results about entry-wise eigenvector analysis for topic models. Section~\ref{sec:TM} applies these results to obtain improved error bounds for the Topic-SCORE algorithm and determine the optimal rate in the short-document case. Section~\ref{sec:proof} describes the main technical components, along with a proof sketch.  Section~\ref{sec:discuss} concludes the paper with discussions. The proofs of all theorems are relegated to the Appendices \ref{appA}--\ref{appE}.

Throughout this paper, 
for a matrix $B$, let $B(i,j)$ or $B_{ij}$ represent the $(i,j)$-th entry. We denote $\|B\|$ as its operator norm and $\|B\|_{2\to \infty}$ as the $2$-to-$\infty$ norm, which is the maximum $\ell_2$ norm across all rows of $B$. For a vector $b$, $b(i)$ or $b_i$ represents the $i$-th component. We denote $\|b\|_1$ and $\|b\|$ as the $\ell_1$ and $\ell_2$ norms of $b$, respectively. The vector ${\bf 1}_n$ stands for an all-one vector of dimension $n$. Unless specified otherwise, $\{e_1, e_2, \ldots, e_p\}$ denotes the standard basis of $\mathbb R^{p}$. Furthermore, we write $a_n \gg b_n$ or $b_n \ll a_n$ if $b_n/a_n = o(1)$ for $a_n, b_n>0$; and we denote $a_n\asymp b_n$ if $C^{-1}b_n< a_n<Cb_n$ for some constant $C>1$.

\section{Entry-Wise Eigenvector Analysis for Topic Models} \label{sec:main}

Let $X\in\mathbb{R}^{p\times n}$ be the word-count matrix following the topic model in \eqref{def:model}. 
We introduce the empirical frequency matrix $D=[d_1,d_2,\ldots,d_n]\in \mathbb{R}^{p\times n}$, defined by:
\begin{equation} \label{def:D}
d_i(j)= N_i^{-1}X_i(j), 
\quad 1\leq i\leq n, 1\leq j \leq p\, .
\end{equation}
Under the model in \eqref{def:model}, we have $\mathbb{E}[d_i]=d_i^0=\sum_{k=1}^Kw_i(k)A_k$. Write $D_0=[d_1^0, d_2^0,\ldots, d_n^0]\in\mathbb{R}^{p\times n}$. It follows that:
\[
\mathbb E D = D_0 = AW.
\]
We observe that $D_0$ is a rank-$K$ matrix; furthermore, the linear space spanned by the first $K$ left singular vectors of $D_0$ is the same as the column space of $A$. \cite{ke2022using} discovered that there is a low-dimensional simplex structure that explicitly connects the first $K$ left singular vectors of $D_0$ with the target topic matrix $A$. This inspired SVD-based methods for estimating $A$.

However, if one directly conducts SVD on $D$, the empirical singular vectors can be noisy because of severe word frequency heterogeneity in natural languages \citep{zipf2013psycho}. In what follows, we first introduce a normalization on $D$ in Section~\ref{sub:model} to handle word frequency heterogeneity and then  derive entry-wise large-deviation bounds for the empirical singular vectors in Section~\ref{sub:entry-wise}. 

\subsection{A Normalized Data Matrix} \label{sub:model}
We first explain why it is inappropriate to conduct SVD on $D$. Let $\bar{N}=n^{-1}\sum_{i=1}^n N_i$ denote the average document length. Write $D=AW+Z$, with $Z=[z_1,z_2,\ldots,z_n] :=D-\mathbb{E}D$. The singular vectors of $D$ are the same as the eigenvectors of $DD'=AWW'A'+AWZ'+ZW'A'+ZZ'$. By model \eqref{def:model}, the columns of $Z$ are centered multinomial-distributed random vectors; moreover, using the covariance matrix formula for multinomial distributions, we have $\mathbb{E}[z_iz_i']=N_i^{-1}[\mathrm{diag}(d_i^0)-d_i^0(d_i^0)']$. It follows that:
\begin{align} \label{normalization-insights}
\mathbb{E}[DD'] & = AWW'A' +  \sum_{i=1}^n N_i^{-1} \bigl[\mathrm{diag}(d_i^0) - d_i^0(d_i^0)'\bigr]\cr
&= AWW'A' + \mathrm{diag}\biggl(\sum_{i=1}^n N_i^{-1}d_i^0 \biggr) - A\biggl(\sum_{i=1}^n N_i^{-1}w_iw_i'\biggr)A'\cr
&= n\cdot A\underbrace{\biggl(\sum_{i=1}^n \frac{N_i-1}{n N_i}w_iw_i'\biggr)}_{\equiv \Sigma_W}A' + \frac{n}{\bar{N}}\cdot \underbrace{\mathrm{diag}\biggl(\sum_{i=1}^n \frac{\bar{N}}{n N_i} d_i^0 \biggr)}_{\equiv M_0}. 
\end{align}
Here $A\Sigma_WA'$ is a rank-$K$ matrix whose eigen-space is the same as the column span of $A$. However, because of the diagonal matrix $M_0$, the eigen-space of $\mathbb{E}[DD']$ is no longer the same as the column span of $A$. 
We notice that the $j$th diagonal of $M_0$ captures the overall frequency of the $j$th word across the whole corpus. Hence, this is an issue caused by word frequency heterogeneity. The second term in \eqref{normalization-insights} is larger when $\bar{N}$ is smaller. This implies that the issue becomes more severe for short documents.

To resolve this issue, we consider a normalization of $D$ to $M_0^{-1/2}D$. It follows that:
\begin{equation} \label{temp:expectation}
\mathbb{E}[M_0^{-1/2}DD'M_0^{-1/2}]= n\cdot M_0^{-1/2}A\Sigma_WA'M_0^{-1/2} + \frac{n}{\bar{N}} I_p. 
\end{equation}
Now, the second term is proportional to an identify matrix and no longer affects the eigen-space. Furthermore, the eigen-space of the first term is the column span of $M_0^{-1/2}A$, and hence, we can use the eigenvectors to recover $M_0^{-1/2}A$ (then $A$ is immediately known).  In practice, $M_0$ is not observed, so we replace it by its empirical version:
 \begin{equation} \label{def:M}
M= {\rm diag}\bigg(\sum_{i=1}^n \frac{\bar{N}}{n N_i} d_i \biggr). 
\end{equation}
We propose to normalize $D$ to $M^{-1/2}D$ before conducting SVD. Later, the singular vectors of $M^{-1/2}D$ will be used in Topic-SCORE to estimate $A$ (see Section~\ref{sec:TM}).

This normalization is similar to the pre-SVD normalization in \cite{ke2022using} but not exactly the same. Inspired by analyzing a special case where $N_i=N$, \cite{ke2022using} proposed to normalize $D$ to $\widetilde{M}^{-1/2}D$, where $\widetilde{M}={\rm diag}(n^{-1}\sum_{i=1}^{n}d_i)$. They continued using $\widetilde{M}$ in general settings, but we discover here that the adjustment of $\widetilde{M}$ to $M$ is necessary when $N_i$'s are unequal.  

\begin{Remark} \label{rmk:trimming}
For extremely low-frequency words, the corresponding diagonal entries of $M$ are very small. This causes an issue when we normalize $D$ to $M^{-1/2}D$. Fortunately, such an issue disappears if we pre-process data. As a standard pre-processing step for topic modeling, we either remove those extremely low-frequency words or combine all of them into a single ``meta-word''. We recommend the latter approach. In detail, let ${\cal L}\subset\{1,2,\ldots,p\}$ be the set of words such that $M(j,j)$ is below a proper threshold $t_n$ (e.g., $t_n$ can be $0.05$ times the average of diagonal entries of $M$). We then sum up all rows of $D$ with indices in ${\cal L}$ to a single row. Let $D^*\in\mathbb{R}^{(p-|{\cal L}|+1)\times n}$ be the processed data matrix. The matrix $D^*$ still has a topic model structure, where each new topic vector results from a similar row combination on the corresponding original topic vector.  
\end{Remark}

\begin{Remark} \label{rmk:Laplacian}
The normalization of $D$ to $M^{-1/2}D$ is reminiscent of the Laplacian normalization in network data analysis, but the motivation is very different. In many network models, the adjacency matrix satisfies that $B=B_0+Y$, where $B_0$ is a low-rank matrix and $Y$ is a generalized Wigner matrix. Since $\mathbb{E}[Y]$ is a zero matrix, the eigen-space of $\mathbb{E}B$ is the same as that of $B_0$. Hence, the role of the Laplacian normalization is not correcting the eigen-space but adjusting the signal-to-noise ratio \citep{ke2022optimal}. In contrast, our normalization here aims to turn $\mathbb{E}[ZZ']$ into an identity matrix (plus a small matrix that can be absorbed into the low-rank part). 
We need such a normalization even under moderate word frequency heterogeneity (i.e., the frequencies of all words are at the same order). 
\end{Remark}

\subsection{Entry-Wise Singular Analysis for $M^{-1/2}D$}\label{sub:entry-wise}
For each $1\leq k\leq K$, let $\hat{\xi}_k\in\mathbb{R}^p$ denote the $k$th left singular vector of $M^{-1/2}D$. Recall that $D_0=\mathbb{E}D$. In addition, define:
\begin{equation} \label{def:M0}
M_0: = \mathbb {E} M  = \mathrm{diag}\biggl(\sum_{i=1}^n \frac{\bar{N}}{n N_i} d_i^0 \biggr). 
\end{equation}
Then, $M_0^{-1/2}D_0$ is a population counterpart of $M^{-1/2}D$. However, the singular vectors of $M_0^{-1/2}D_0$ are not the population counterpart of $\hat{\xi}_k$'s. In light of \eqref{temp:expectation}, we define:
\begin{equation} \label{def:popEigVec}
\xi_k: \mbox{the $k$th eigenvector of } M_0^{-1/2}\mathbb{E}[DD']M_0^{-1/2}, \qquad 1\leq k\leq K. 
\end{equation}
Write $\hat \Xi : = [\hat \xi_1, \cdots, \hat \xi_K]$ and $\Xi: = [\xi_1, \cdots,  \xi_K]$. We aim to derive a large-deviation bound for each individual row of $(\hat{\Xi}-\Xi)$, subject to a column rotation of $\hat{\Xi}$.

We need a few assumptions. Let $h_j=\sum_{k=1}^K A_k(j)$ for $1\leq j\leq p$. Define:
\begin{equation} \label{notation}
H = {\rm diag} (h_1, \cdots, h_p), \qquad \Sigma_A = A'H^{-1}A, \qquad \Sigma_W = \frac{1}{n}\sum_{i=1}^n (1-N_i^{-1}) w_iw_i'. 
\end{equation}
Here $\Sigma_A$ and $\Sigma_W$ are called the topic-topic overlapping matrix and the topic-topic concurrence matrix, 
respectively, \citep{ke2022using}. It is easy to see that $\Sigma_W$ is properly scaled. We remark that $\Sigma_A$ is also properly scaled, because $\sum_{\ell=1}^K \Sigma_A(k,\ell)= \sum_{j=1}^p \sum_{\ell=1}^K h_j^{-1}A_k(j)A_\ell(j)=1$. 
\begin{asm} \label{asm:h}
Let $h_{\max} = \max_{1\leq j\leq p} h_j $, $h_{\min} = \min_{1\leq j\leq p } h_j$ and $\bar{h}= \frac 1p \sum_{j=1}^p h_j$. We assume:
\[
h_{\min} \geq c_1 \bar{h} = c_1K/p, \qquad \text{ for a constant $c_1\in (0,1)$}. 
\]
\end{asm} 
\begin{asm}\label{asm:WA}
For a constant $c_2\in (0,1)$ and a sequence $\beta_n\in (0,1)$, we assume:
\[
\lambda_{\min} (\Sigma_W) \geq c_2, \qquad \lambda_{\min} (\Sigma_A)\geq c_2\beta_n, \qquad \min_{1\leq k,\ell \leq K} \Sigma_A(k, \ell) \geq c_2.
\]
\end{asm}

Assumption~\ref{asm:h} is related to word frequency heterogeneity. Each $h_j$ captures the overall frequency of word $j$, and $\bar{h}=p^{-1}\sum_j h_j = p^{-1}\sum_k \|A_k\|_1 = K/p$. 
By Remark~\ref{rmk:trimming}, all extremely low-frequency words have been combined in pre-processing. It is reasonable to assume that $h_{\min}$ is at the same order of $\bar{h}$. Meanwhile, we put no restrictions here on $h_{\max}$, so that $h_j$'s can still be at different orders.

Assumption~\ref{asm:WA} is about topic weight balance and between-topic similarity. $\Sigma_W$ can be regarded as an affinity matrix of $w_i$'s. It is mild to assume that $\Sigma_W$ is well-conditioned. In a special case where $N_i=N$ and each $w_i$ is degenerate, $\Sigma_W$ is a diagonal matrix whose $k$th diagonal entry is the fraction of documents that put all weights on topic $k$; hence, $\lambda_{\min}(\Sigma_W)\geq c_2$ is interpreted as ``topic weight balance''.  
Regarding $\Sigma_A$, we have seen that it is properly scaled (its maximum eigenvalue is at the constant order). When $K$ topic vectors are exactly the same, $\lambda_{\min}(\Sigma_A)=0$; when the topic vectors are not the same, 
$\lambda_{\min}(\Sigma_A)\neq 0$, and it measures the signal strength.  \cite{ke2022using} assumed that $\lambda_{\min}(\Sigma_A)$ is bounded below by a constant, but we allow weaker signals by allowing $\lambda_{\min}(\Sigma_A)$ to diminish as $n\to\infty$. 
We also require a lower bound on $\Sigma_A(k, \ell) $, meaning that there should be certain overlaps between any two topics. This is reasonable as some commonly used words are not exclusive to any one topic and tend to occur frequently \citep{ke2022using}.

The last assumption is about the vocabulary size and document lengths.  
\begin{asm}\label{asm:para}
There exists $N\geq 1$ and a constant $c_3\in (0,1)$ such that $c_3 N\leq N_i\leq c_3^{-1}N$ for all $1\leq i\leq n$. In addition, for an arbitrary constant $C_0>0$:
\[
\min\{p, N\} \geq \log^3 (n),   \qquad \max\{\log(p), \log( N)\} \leq {C_0}\log(n), \qquad p \log^{2} (n) {\leq}  Nn\beta^2_n.
\]
\end {asm}

In Assumption \ref{asm:para}, the first two inequalities restrict that $N$ and $p$ are between $\log^3(n)$ and $n^{C_0}$, for an arbitrary constant $C_0>0$. This covers a wide regime, including the scenarios of both long documents ($N\geq c\cdot p$) and short documents ($N=o(p)$). 
The third inequality is needed so that the canonical angles between the empirical and population singular spaces converge to zero, which is necessary for our singular vector analysis. This condition is mild, as $Nn$ is the order of total word count in the corpus, which is often much larger than $p$.

With these assumptions, we now present our main theorem.
\begin{Theorem}[Entry-wise singular vector analysis]\label{thm:row_eigenv}
Fix $K\geq 2$ and positive constants $c_1, c_2, c_3$, { and $C_0$}. 
Under the model (\ref{def:model}), suppose Assumptions \ref{asm:h}--\ref{asm:para} hold. {For any constant $C_1>0$, there exists $C_2>0$ such that} with probability $1- {n^{-C_1}}$, there is an orthogonal matrix $O\in \mathbb{R}^{K\times K}$ satisfying that simultaneously for $1\leq j \leq p$:
\[
\Vert e_j' (\hat \Xi- \Xi O')   \Vert \leq {C_2} \sqrt{\frac{h_j p\log (n)}{nN\beta_n^2}}. 
\]
The constant $C_2$ only depends on $C_1$ and $(K, c_1,c_2,c_3, C_0)$.
\end{Theorem}
In Theorem~\ref{thm:row_eigenv}, we do not assume any gap among the $K$ singular values of $M_0^{-1/2}D_0$; hence, it is only possible to recover $\Xi$ up to a column rotation $O$. 
The sin-theta theorem \citep{sin-theta} enables us to bound $\|\hat{\Xi}-\Xi O'\|^2_F=\sum_{j=1}^p \|e_j'(\hat{\Xi}-\Xi O')\|^2$, but it is insufficient for analyzing spectral algorithms for topic modeling (see Section~\ref{sec:TM}). We need a bound for each individual row of $(\hat{\Xi}-\Xi O')$, and this bound should depend on $h_j$ properly.  

We compare Theorem~\ref{thm:row_eigenv} with the result in \cite{ke2022using}. They assumed that $\beta_n^{-1}=O(1)$, so their results are only for the strong-signal regime. They showed that when $n$ is sufficiently large:
\begin{equation} \label{original-rate}
\Vert e_j' (\hat \Xi- \Xi O')   \Vert \leq C\biggl(1+\min\Bigl\{\frac{p}{N},\ \frac{p^2}{N\sqrt{N}}\Bigr\}\biggr)\sqrt{\frac{h_j p\log (n)}{nN}}.
\end{equation}

When $N\geq c\cdot p$ (long documents), it is the same bound as in Theorem~\ref{thm:row_eigenv} (with $\beta_n=1$). However, when $N=o(p)$ (short documents), it is strictly worse than Theorem~\ref{thm:row_eigenv}. 
We obtain better bounds than those in \cite{ke2022using} because of new proof ideas, especially the use of refined perturbation analysis for SVD and a decoupling technique for U-statistics (see Section~\ref{sub:technical}). 



\section{Improved Rates for Topic Modeling} \label{sec:TM}
We apply the results in Section~\ref{sec:main} to improve the error rates of topic modeling.  
Topic-SCORE  \citep{ke2022using} is a spectral algorithm for estimating the topic matrix $A$. It achieves the optimal rate in the long-document case ($N\geq c\cdot p$). However, in the short-document case ($N=o(p)$), the known rate of Topic-SCORE does not match with the minimax lower bound. We address this gap by providing better error bounds for Topic-SCORE. Our results reveal the optimal rate for topic modeling in the short-document case for the first time.  


\subsection{The Topic-Score Algorithm} \label{sub:T-SCORE}
Let $\hat{\xi}_1,\hat{\xi}_2,\ldots,\hat{\xi}_K$ be as in Section~\ref{sec:main}. 
Topic-SCORE first obtains word embeddings from these singular vectors. 
Note that $M^{-1/2}D$ is a non-negative matrix. By Perron's theorem~\citep{HornJohnson}, under mild conditions, $\hat{\xi}_1$ is a strictly positive vector. Define $\hat{R}\in\mathbb{R}^{p\times (K-1)}$ by:
\begin{equation} \label{SCORE}
\hat{R}(j,k) = \hat{\xi}_{k+1}(j)/\hat{\xi}_1(j), \qquad 1\leq j\leq p, 1\leq k\leq K-1. 
\end{equation}
Let $\hat{r}_1',\hat{r}_2',\ldots,\hat{r}_p'$ denote the rows of $\hat{R}$. Then, $\hat{r}_j$ is a $(K-1)$-dimensional embedding of the $j$th word in the vocabulary. This is known as the SCORE embedding \cite{SCORE, SCOREreview}, which is now widely used in analyzing heterogeneous network and text data. 

\begin{figure}[tb!]
\centering
\includegraphics[width=.8\textwidth]{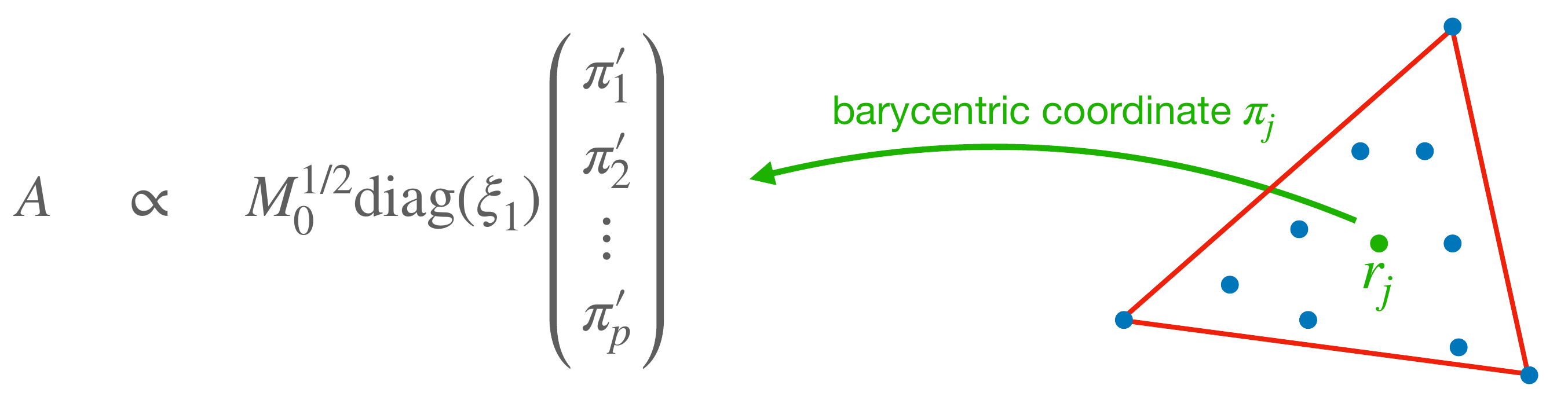}
\caption{An illustration of Topic-SCORE in the noiseless case ($K=3$). The blue dots are $r_j\in\mathbb{R}^{K-1}$ (word embeddings  constructed from the population singular vectors), and they are contained in a simplex with $K$ vertices. This simplex can be recovered from a vertex hunting algorithm. Given this simplex, each $r_j$ has a unique barycentric coordinate $\pi_j\in\mathbb{R}^K$. The topic matrix $A$ is recovered from stacking together these $\pi_j$'s and utilizing $M_0$ and $\xi_1$.} \label{fig:T-SCORE}
\end{figure}

\cite{ke2022using} discovered that there is a simplex structure associated with these word embeddings. 
Specifically, let $\xi_1,\xi_2,\ldots,\xi_K$ be the same as in \eqref{def:popEigVec} and define the population counterpart of $\hat{R}$ as $R$, where:
\begin{equation} \label{SCORE-oracle}
R(j,k) = \xi_{k+1}(j)/\xi_1(j), \qquad 1\leq j\leq p, 1\leq k\leq K-1. 
\end{equation}
Let $r_1', r_2',\ldots, r_p'$ denote the rows of $R$. All these $r_j$ are contained in a simplex ${\cal S}\subset\mathbb{R}^{K-1}$ that has $K$ vertices $v_1, v_2, \ldots, v_K$ (see Figure~\ref{fig:T-SCORE}).  
If the $j$th word is an anchor word \citep{Ge,donoho2003does} (an anchor word of topic $k$ satisfies that $A_k(j)\neq 0$ and $A_\ell(j)=0$ for all other $\ell\neq k$), then $r_j$ is located at one of the vertices. Therefore, as long as each topic has at least one anchor word, we can apply a vertex hunting \citep{ke2022using} algorithm to recover the $K$ vertices of ${\cal S}$. By definition of a simplex, each point inside ${\cal S}$ can be written uniquely as a convex combination of the $K$ vertices, and the $K$-dimensional vector consisting of the convex combination coefficients is called the barycentric coordinate.
After recovering the vertices of ${\cal S}$, we can easily compute the barycentric coordinate $\pi_j\in\mathbb{R}^K$ for each $r_j$. Write $\Pi=[\pi_1,\pi_2,\ldots,\pi_p]'$. \cite{ke2022using} showed that:
\[
A_k \;\;\propto\;\; M_0^{1/2}\mathrm{diag}(\xi_1)\Pi e_k, \qquad 1\leq k\leq K. 
\]
Therefore, we can recover $A_k$ by taking the $k$th column of $M_0^{1/2}\mathrm{diag}(\xi_1)\Pi$ and re-normalizing it to have a unit $\ell^1$-norm. This gives the main idea behind Topic-SCORE (see Figure~\ref{fig:T-SCORE}). 

The full algorithm is given in Algorithm~\ref{alg:T-SCORE}. It requires plugging in a vertex hunting (VH) algorithm. 
A VH algorithm aims to estimate $v_1,v_2,\ldots, v_K$ from the noisy point cloud $\{\hat{r}_j\}_{1\leq j\leq p}$.  There are many existing VH algorithms (see sec 3.4 of \cite{SCOREreview}). A VH algorithm is said to be efficient if it satisfies that $\max_{1\leq k\leq K}\|\hat{v}_k-v_k\|\leq C\max_{1\leq j\leq p}\|\hat{r}_j-r_j\|$ (subject to a permutation of $\hat{v}_1,\hat{v}_2,\ldots,\hat{v}_K$). 
We always plug in an efficient VH algorithm, such as the successive projection algorithm  \citep{araujo2001successive}, the pp-SPA algorithm \citep{SP-ICLR}, and several algorithms in sec 3.4 of \cite{SCOREreview}.

\begin{algorithm}[htb]
\caption{Topic-SCORE} \label{alg:T-SCORE}
{\bf Input}: D,  K, and a vertex hunting (VH) algorithm. 
\begin{itemize}
\item (\textit{Word embedding})  Let $M$ be as in \eqref{def:M}. Obtain $\hat \xi_1, \hat\xi_2, \cdots, \hat \xi_K$, the first $K$ left singular vectors of $M^{-1/2}D$. Compute $\hat R$ as in \eqref{SCORE} and write $\hat R=[\hat r_1, \hat r_2, \cdots, \hat r_p]'$. 
\item (\textit{Vertex hunting}). Apply the VH algorithm on $\{\hat r_j\}_{1\leq j\leq p}$ to get $\hat v_1, \cdots, \hat v_K$. 
\item (\textit{Topic matrix estimation}) For $1\leq j \leq p$, solve $\hat \pi_j^*$ from:
\[ 
\left(\begin{array} {ccc}
1&\ldots & 1\\
\hat v_1 & \ldots & \hat v_K
\end{array}
\right)\hat \pi_j^* = \left ( \begin{array} {c} 1\\\hat r_j\end{array}\right). 
\]
Let $\tilde \pi^*_j =\max\{\hat \pi_j^*,  0 \}$ (the maximum is taken component-wise) and $\hat \pi_i = \tilde \pi^*_j/ \|\tilde \pi^*_j\|_1 $. Write $\hat \Pi = [\hat \pi_1, \ldots, \hat\pi_p]'$. Let  $\tilde A = M^{1/2} {\rm diag}(\hat \xi_1) \hat \Pi$. Obtain $\hat A =\tilde A [{\rm diag} ({\bf 1}_p' \tilde A) ]^{-1}$. 
\end{itemize}
{\bf Output}: the estimated topic matrix $ \hat A$. 
\end{algorithm}

Additionally, after $\hat A$ is obtained, \cite{ke2022using} suggested to estimate $w_1,w_2,\ldots,w_n$ as follows. We first run a weighted least-squares to obtain $\hat{w}_i^*$: 
\begin{equation}\label{estimateW}
\hat{w}^*_i = \mathrm{argmin}_{w\in\mathbb{R}^K}\|M^{-1/2}(d_i - Aw_i) \|^2, \qquad 1\leq i\leq n. 
\end{equation}
Then, set all the negative entries of $\hat{w}_i^*$ to zero and re-normalize the vector to have a unit $\ell^1$-norm. The resulting vector is $\hat{w}_i$. 

{

\begin{Remark}
In real-world applications, both $n$ and $p$ can be very large. However, since $\hat{R}$ is constructed from only a few singular vectors, its rows are only in dimension $(K-1)$. It suggests that Topic-SCORE leverages a `low-dimensional' simplex structure and is scalable to large datasets. When $K$ is bounded, the complexity of Topic-SCORE is at most $O(np\min\{n,p\})$ \citep{ke2022using}. 
The real computing time was also reported in \cite{ke2022using} for various values of ($n,p$). For example, when both $n$ and $p$ are a few thousands, it takes only a few seconds to run Topic-SCORE.  
\end{Remark}

}

\subsection{The Improved Rates for Estimating $A$ and $W$}\label{sub:rates-topic}

We provide the error rates of Topic-SCORE. First, we study the word embeddings  $\hat{r}_j$.  By \eqref{SCORE}, $\hat{r}_j$ is constructed from the $j$th row of $\hat{\Xi}$. Therefore, we can apply Theorem~\ref{thm:row_eigenv} to derive a large-deviation bound for $\hat{r}_j$. 

Without loss of generality, we set $C_1=4$ henceforth, transforming the event probability $1-n^{-C_1}$ in Theorem~\ref{thm:row_eigenv} to $1-o(n^{-3})$. 
We also use $C$ to denote a generic constant, whose meaning may change from one occurrence to another. In all instances, $C$ depends sorely on $K$ and the constants $(c_1, c_2, c_3, C_0)$ in Assumptions \ref{asm:h}--\ref{asm:para}. 
\begin{Theorem}[Word embeddings]\label{thm:row_R}
Under the setting of Theorem~\ref{thm:row_eigenv}, with probability $1-o(n^{-3})$, 
there exists an orthogonal matrix $\Omega \in \mathbb{R}^{(K-1)\times (K-1)}$ such that simultaneously for $1\leq j \leq p$:
\[
\Vert \hat r_j - \Omega r_j   \Vert \leq C \sqrt{\frac{p\log(n)}{nN\beta_n^2}}. 
\]
\end{Theorem}

Next, we study the error of $\hat{A}$. The $\ell^1$-estimation error is $
{\cal L}(\hat{A}, A):=\sum_{k=1}^{K}\|\hat{A}_k-A_k\|_1$, subject to an arbitrary column permutation of $\hat{A}$.
For ease of notation, we do not explicitly denote this permutation in theorem statements, but it is accounted for in the proofs.   For each $1\leq j\leq p$, let $\hat{a}_j'\in\mathbb{R}^K$ and $a_j'\in\mathbb{R}^K$ denote the $j$th row of $\hat A$ and $A$, respectively. We can re-write the $\ell^1$-estimation error as ${\cal L}(\hat{A}, A)=\sum_{j=1}^{p}\|\hat{a}_j-a_j\|_1$. The next theorem provides an error bound for each individual $\hat a_j$, and the aggregation of these bounds yields an overall bound for ${\cal L}(\hat{A}, A)$:

\begin{Theorem}[Estimation of $A$]\label{thm:A_k}
Under the setting of Theorem~\ref{thm:row_eigenv}, we additionally assume that each topic has at least one anchor word. 
With probability $1- o(n^{-3})$,  simultaneously for $1\leq j \leq p$:
\[
\Vert \hat a_j - a_j  \Vert_1 \leq \|a_j\|_1\cdot C \sqrt{\frac{p\log (n)}{nN\beta_n^2}}\,. 
\]
Furthermore, with probability $1- o(n^{-3})$, the $\ell^1$-estimation error satisfies that:
\[
\mathcal L (\hat A, A)  
\leq C  \sqrt{\frac{p\log (n)}{nN\beta_n^2}}.
\]
\end{Theorem}

Theorem~\ref{thm:A_k} improves the result in \cite{ke2022using} in two aspects. First, \cite{ke2022using} assumed $\beta_n^{-1}=O(1)$, so their results did not allow for weak signals. Second, even when $\beta_n^{-1}=O(1)$, their bound is worse than ours by a factor similar to that in \eqref{original-rate}.

Finally, we have the error bound for estimating $w_i$'s using the estimator in \eqref{estimateW}. 
\begin{Theorem}[Estimation of $W$]\label{thm:w_i}
Under the setting of Theorem~\ref{thm:A_k}, with probability $1- o(n^{-3})$, 
subject to a column permutation of $\hat{W}$:
\[
\max_{1\leq i\leq n}\Vert \hat w_i-  w_i   \Vert_1 \leq C \beta_n^{-1} \Biggl(\sqrt{\frac{p\log (n)}{nN\beta_n^2}} + C\sqrt{\frac{\log(n)}{N}}\Biggr).  
\]
\end{Theorem}

In Theorem~\ref{thm:w_i}, there are two terms in the error bound of $\hat{w}_i$. The first term comes from the estimation error in $\hat A$, and the second term is from noise in $d_i$. In the strong-signal case of $\beta_n^{-1}=O(1)$, we can compare Theorem~\ref{thm:w_i} with the bound for $\hat{w}_i$ in  \cite{ke2022using}. The bound there also has two terms: its second term is similar to ours, but its first term is strictly worse.

\subsection{Connections and Comparisons} \label{sub:discuss}

There have been numerous results about the error rates of estimating $A$ and $W$. For example,  \citep{Ge} provided the first explicit theoretical guarantees for topic modeling, but they did not study the statistical optimality of their method. Recently, the statistical literature aimed to understand the fundamental limits of topic modeling. Assuming $\beta_n^{-1}=O(1)$, refs.~\citep{ke2022using,bing2020fast} gave a minimax lower bound, $\sqrt{p/(Nn)}$, for the rate of estimating $A$, and refs.~\citep{wu2023sparse,klopp2023assigning} gave a minimax lower bound, $1/\sqrt{N}$, for estimating each $w_i$.

For estimating $A$, when $\beta_n^{-1}=O(1)$, the existing theoretical results are summarized in Table~\ref{tb:rates}. 
Ours is the only one that matches the minimax lower bound across the entire regime.   
In the long-document case ($N\geq c\cdot p$, Cases 1--2 in Table~\ref{tb:rates}), the error rates in \cite{ke2022using,bing2020fast} together have matched the lower bound, concluding that $\sqrt{p/(Nn)}$ is indeed the optimal rate. However, in the short-document case ($N=o(p)$, Case 3 in Table~\ref{tb:rates}), there was a gap between the lower bound and the existing error rates.  
Our result addresses the gap and concludes that $\sqrt{p/(Nn)}$ is still the optimal rate. 
When $\beta_n=o(1)$, the error rates of estimating $A$ were rarely studied. We conjecture that $\sqrt{p/(Nn\beta_n^2)}$ is the optimal rate, and the Topic-SCORE algorithm is still rate-optimal. 

\begin{table}[tb!]
\centering
\caption{A summary of the existing theoretical results for estimating $A$ ($n$ is the number of documents, $p$ is the vocabulary size, $N$ is the order of document lengths, and $h_{\max}$ and $h_{\min}$ are the same as in \eqref{notation}). Cases 1--3 refer to $N\geq p^{4/3}$, $p\leq N<p^{4/3}$, and $N<p$, respectively. For Cases 2--3, the sub-cases `a' and `b' correspond to $n\geq \max\{ Np^2, p^3, N^2 p^5\}$ and $n< \max\{ Np^2, p^3, N^2 p^5\}$, respectively. We have translated the results in each paper to the bounds on ${\cal L}(\hat A, A)$, with any logarithmic factor omitted. } \label{tb:rates}
\scalebox{.9}{
\begin{tabularx}{\textwidth}{lccccc} 
\toprule
&\textbf{\boldmath{Case 1}}& \textbf{\boldmath{Case 2a}}& \textbf{\boldmath{Case 2b}}& \textbf{\boldmath{Case 3a}}& \textbf{\boldmath{Case 3b}} \\ 
\midrule
\multirow{1.3}{*}{\cite{ke2022using}}&  $\sqrt{\frac{p}{Nn}}$ &  $\sqrt{\frac{p}{Nn}}$ &  $\frac{p^2}{N\sqrt N}\sqrt{\frac{p}{Nn}}$ &   $\frac{p}{N}\sqrt{\frac{p}{Nn}} $ & $ \frac{p^2}{N\sqrt N}\sqrt{\frac{p}{Nn}}$\\ 
\midrule
\multirow{1.3}{*}{\cite{Ge}}& $\frac{p^4}{\sqrt{Nn}}$ &  $\frac{p^4}{\sqrt{Nn}}$ &  $\frac{p^4}{\sqrt{Nn}}$ &   $\frac{p^4}{\sqrt{Nn}}$ & $ \frac{p^4}{\sqrt{Nn}}$\\ 
\midrule 
\multirow{1.3}{*}{\cite{bing2020fast}}& $\sqrt{\frac{p}{Nn}} \cdot \frac{h_{\max}}{h_{\min}} $ &  $\sqrt{\frac{p}{Nn}}\cdot\frac{h_{\max}}{h_{\min}} $ & $\sqrt{\frac{p}{Nn}} \cdot\frac{h_{\max}}{h_{\min}} $ & \multirow{1.3}{*}{{NA} 
}
 & \multirow{1.3}{*}{NA}\\
\midrule
\multirow{1.3}{*}{\cite{bansal2014provable}}& $N\sqrt{\frac pn}$ &  $N\sqrt{\frac pn}$ &  $N\sqrt{ \frac pn}$ &   $N\sqrt{ \frac pn} $ & $N\sqrt{\frac pn}$ \\ 
\midrule
\multirow{1.5}{*}{Our results} & $\sqrt{\frac{p}{Nn}}$ &  $\sqrt{\frac{p}{Nn}}$ &  $\sqrt{\frac{p}{Nn}}$ &   $\sqrt{\frac{p}{Nn}} $ & $\sqrt{\frac{p}{Nn}}$ \\ [0.2cm]
\bottomrule
\end{tabularx}}
\end{table}

We emphasize that our rate is not affected by severe word frequency heterogeneity. As long as $h_{\min}/\bar{h}$ is lower bounded by a constant (see Assumption~\ref{asm:h} and explanations therein), our rate is always the same, regardless of the magnitude of $h_{\max}$. 
In contrast, the error rate in \cite{bing2020fast} is sensitive to word frequency heterogeneity, with an extra factor of $h_{\max}/h_{\min}$ that can be as large as $p$. There are two reasons that enable us to achieve a flat rate even under severe word frequency heterogeneity: one is the proper normalization of data matrix, as described in Section~\ref{sub:model}, and the other is the careful analysis of empirical singular vectors (see Section~\ref{sec:proof}).

For estimating $W$, when $\beta_n^{-1}=O(1)$, our error rate in Theorem~\ref{thm:w_i} matches the minimax lower bound if $n\geq p\log(n)$. Our approach to estimating $W$ involves first obtaining $\hat A$ and then regressing $d_i$ on $\hat A$ to derive $\hat{w}_i$. The condition $n\geq p\log(n)$ ensures that the estimation error in $\hat{A}$ does not dominate the overall error. This condition is often met in scenarios where a large number of documents can be collected, but the vocabulary size remains relatively stable. However, if $n<p\log(n)$, a different approach is necessary, requiring the estimation of $W$ first. This involves using the right singular vectors of $M^{-1/2}D$. While our analysis has primarily focused on the left singular vectors, it can be extended to study the right singular vectors as well.


\section{Proof Ideas}\label{sec:proof}
Our main result is Theorem~\ref{thm:row_eigenv}, which provides entry-wise large-deviation bounds for singular vectors of $M^{-1/2}D$. Given this theorem, the proofs of Theorems~ \ref{thm:row_R}--\ref{thm:w_i} are similar to those in \citep{ke2022using} and thus relegated to the appendix. In this section, we focus on discussing the proof techniques of Theorem~\ref{thm:row_eigenv}. 

\subsection{Why the Leave-One-Out Technique Fails}
Leave-one-out \citep{abbe2020entrywise,ke2022optimal} is a common technique in entry-wise eigenvector analysis for a Wigner-type random matrix $B=B_0+Y\in\mathbb{R}^{m\times m}$, where $B_0$ is a  symmetric non-stochastic low-rank matrix and $Y$ is a symmetric random matrix whose upper triangle consists of independent mean-zero variables. 
One example of such matrices is the adjacency matrix of a random graph generated from the block-model family \citep{SCORE}. 

However, our target here is the singular vectors of $M^{-1/2}D$, which are the eigenvectors of $B:=M^{-1/2}DD'M^{-1/2}$. 
This is a Wishart-type random matrix, whose upper triangular entries are not independent. 
We may also construct  a symmetric matrix:
\[
\mathcal G: = \begin{pmatrix}
0  & M^{-1/2}D\\
D'M^{-1/2} & 0
\end{pmatrix} \quad \in\quad\mathbb{R}^{(p+n)\times (p+n)}. 
\]

The eigenvectors of ${\cal G}$ take the form $\hat u_k = (\hat \xi_k', \hat \eta_k')'$, $1\leq k\leq K$, where $\hat{\xi}_k\in\mathbb{R}^p$ and $\hat\eta_k\in\mathbb{R}^n$ are the $k$th left and right singular vectors of $M^{-1/2} D$, respectively. Unfortunately, the upper triangle of ${\cal G}$ still contains dependent entries. 
Some dependence is from the normalization matrix $M$. It may be addressed by using the techniques developed by \cite{ke2022optimal} in studying graph Laplacian matrices. 
A more severe issue is the dependence among entries in $D$. According to basic properties of multinomial distributions, 
$D$ only has column independence but no row independence. As a result, even after we replace $M$ by $M_0$, the $j$th row and column of ${\cal G}$ are still dependent of the remaining ones, for each $1\leq j\leq p$. In conclusion, we cannot apply the leave-one-out technique on either $B$ or ${\cal G}$.


\subsection{The Proof Structure in \cite{ke2022using} and Why It Is Not Sharp for Short Documents} \label{sub:technical}

Our entry-wise eigenvector analysis primarily follows the proof structure in \cite{ke2022using}. Recall that $\hat{\xi}_k\in\mathbb{R}^p$ is the $k$th left singular vector of $M^{-1/2}D$.  
Define:
\begin{equation}\label{def:GG0}
G : = M^{-1/2} DD' M^{-1/2} - \frac{n}{\bar{N}} I_p, \qquad G_0:=  n\cdot M_0^{-1/2} A\Sigma_W A'M_0^{-1/2}\,. 
\end{equation}

Since the identify matrix in $G$ does not affect the eigenvectors, $\hat{\xi}_k$ is the $k$th eigenvector of $G$. Additionally, 
it follows from \eqref{def:popEigVec} and \eqref{temp:expectation} that $\xi_k$ is the $k$th eigenvector of $G_0$. 
By \eqref{temp:expectation}:
\begin{equation} \label{temp:G-G0}
G - G_0 = M^{-1/2}DD'M^{-1/2}-M_0^{-1/2}\mathbb{E}[DD']M_0^{-1/2}. 
\end{equation}

The entry-wise eigenvector analysis in \cite{ke2022using} has two steps. 
\begin{itemize}
\item Step 1: Non-stochastic perturbation analysis. In this step, no distributional assumptions are made on $G$. The analysis solely focuses on connecting the perturbation from $\Xi$ to $\hat{\Xi}$ with the perturbation from $G_0$ to $G$. 
They showed in Lemma F.1 \citep{ke2022using}: 
 \begin{equation} \label{temp:OldProof-1}
\| e_j'( \hat \Xi - \Xi O' ) \|\leq C\|G_0\|^{-1} \big( \| e_j'\Xi\| \|G- G_0\| + \sqrt K \| e_j' (G- G_0)\|\big). 
\end{equation}
\item Step 2: Large-deviation analysis of $G-G_0$. In this step, \cite{ke2022using} derived the large-deviation bounds for $\|G-G_0\|$ and $\|e_j'(G-G_0)\|$ under the multinomial model \eqref{def:model}. 
For example, they showed in Lemma F.5 \citep{ke2022using} that when $n$ is properly large, with high probability:
\begin{equation} \label{temp:OldProof-2}
\|G-G_0\|\leq C\bigl(1+N^{-1}\sqrt{p}\bigr)\sqrt{\frac{np\log(n)}{N}}. 
\end{equation} 
\end{itemize}

However, when $N=o(p)$ (short documents), neither step is sharp.  
In  \eqref{temp:OldProof-1}, the second term $\|e_j'(G-G_0)\|$ was introduced as an upper bound for $\|e_j' (G- G_0)\hat \Xi\|$, but this bound is too crude. In Section~\ref{subsec:perturbation}, we will conduct careful analysis of $\|e_j' (G- G_0)\hat \Xi\|$ and introduce a new perturbation bound which significantly improves \eqref{temp:OldProof-1}. 
In \eqref{temp:OldProof-2}, the spectral norm is controlled via an $\varepsilon$-net argument \citep{Vershynin}, which reduces the analysis to studying a quadratic form of $Z$; \cite{ke2022using} analyzed this quadratic form by applying martingale Bernstein inequality. Unfortunately, in the short-document case, it is hard to control  the conditional variance process of the underlying martingale. In Section~\ref{subsec:U-statistics}, we address it by leveraging the matrix Bernstein inequality \citep{tropp2012user} and the decoupling inequality \citep{de1995decoupling, de2012decoupling} for U-statistics.

\subsection{Non-Stochastic Perturbation Analysis} \label{subsec:perturbation}
In this subsection, we abuse notations to use $G$ and $G_0$ to denote two arbitrary $p\times p$ symmetric matrices, with $\mathrm{rank}(G_0)=K$.  For $1\leq k\leq K$, let $\hat{\lambda}_k$ and $\lambda_k$ be the $k$th largest eigenvalue (in magnitude) of $G$ and $G_0$, respectively, and let $\hat{\xi}_k\in\mathbb{R}^p$ and $\xi_k\in\mathbb{R}^p$ be the associated eigenvectors. Write $\hat \Lambda  = {\rm diag} (\hat \lambda_1, \hat\lambda_2, \ldots, \hat \lambda_K)$, $\hat\Xi= [\hat{\xi}_1, \hat\xi_2, \ldots,\hat{\xi}_K]$, and define $\Lambda$ and $\Xi$ similarly. 
Let $U\in\mathbb{R}^{K\times K}$  and $V\in\mathbb{R}^{K\times K}$  be such that its columns contain the left and right singular vectors of $\hat\Xi'\Xi$, respectively. Define ${\rm sgn} (\hat \Xi' \Xi) = U'V$. 
For any matrix $B$ and $q>0$, let $\|B\|_{q\to\infty}=\max_{i}\|e_i'B\|_q$.  

\begin{Lemma}\label{lem:perturb}
Suppose $\|G-G_0\|\leq (1-c_0)|\hat{\lambda}_K|$, for some $c_0\in (0,1)$. 
Consider an arbitrary $p\times p$ diagonal matrix $\Gamma={\rm diag}(\gamma_1,\gamma_2, \ldots,\gamma_p)$, where:
\[
\gamma_j>0 \mbox{ is an upper bound for } \|e_j'\Xi\|\|G-G_0\| + \|e_j'(G-G_0)\Xi\|. 
\]
If $\|\Gamma^{-1}(G-G_0)\Gamma\|_{1\to\infty}\leq (1-c_0)|\hat{\lambda}_K|$, then for the orthogonal matrix $O = {\rm sgn} (\hat \Xi' \Xi)$, it holds simultaneously for $1\leq j\leq p$ that:
\[
\|e_j'(\hat{\Xi} -\Xi O')\|\leq c_0^{-1}|\hat{\lambda}_K|^{-1}\gamma_j.
\] 
\end{Lemma}

Since $\gamma_j$ is an upper bound for $\|e_j'\Xi\|\|G-G_0\| + \|e_j'(G-G_0)\Xi\|$, we can interpret the result in  Lemma~\ref{lem:perturb} as:
\begin{equation} \label{main:perturb}
\|e_j'(\hat{\Xi} -\Xi O')\|  \leq C|\hat{\lambda}_K|^{-1}\bigl(\|e_j'\Xi\|\|G-G_0\| + \|e_j'(G-G_0)\Xi\|\bigr). 
\end{equation} 

Comparing \eqref{main:perturb} with \eqref{temp:OldProof-1}, the second term has been reduced. Since $\Xi$ projects the vector $e_j'(G-G_0)$ into a much lower dimension, we expect that $\|e_j'(G-G_0)\Xi\|\ll \|e_j'(G-G_0)\|$ in many random models for $G$. In particular, this is true for the $G$ and $G_0$ defined in \eqref{def:GG0}. Hence, there is a significant improvement over the analysis in \cite{ke2022using}.


\subsection{Large-Deviation Analysis of $(G-G_0)$} \label{subsec:U-statistics}

In this subsection, we focus on the specific $G$ and $G_0$ as defined in \eqref{def:GG0}. The crux of proving Theorem~\ref{thm:row_eigenv} lies in determining the upper bound $\gamma_j$ as defined in Lemma~\ref{lem:perturb}. This is accomplished through the following lemma.

\begin{Lemma} \label{lem:large-dev-G}
Under the settings of Theorem \ref{thm:row_eigenv}, let $G$ and $G_0$ be as in (\ref{def:GG0}). 
For any constant $C_1>0$, there exists $C_3>0$ such that with probability $1- n^{-C_1}$, 
simultaneously for $1\leq j\leq p$:
\[
\Vert G- G_0 \Vert  \leq  C_3  \sqrt{\frac{pn\log(n)}{N}}, \qquad \Vert e_j'(G-G_0)\Xi\Vert\leq  C_3  \sqrt{\frac{h_j np\log(n)}{N}}. 
\] 
The constant $C_3$  only depends on $C_1$ and $(K, c_1, c_2, c_3, C_0)$.
\end{Lemma}

We compare the bound for $\|G-G_0\|$ in Lemma~\ref{lem:large-dev-G} with the one in~\cite{ke2022using} as summarized in~\eqref{temp:OldProof-2}. There is a significant improvement when $N\leq p^2$. 
This improvement primarily stems from the application of a decoupling inequality for U-statistics, as elaborated below.

We outline the proof of the bound for $\|G-G_0\|$. Let $Z=D-\mathbb{E}[D]=[z_1,z_2,\ldots,z_n]$. From \eqref{eq:G-G_0}-\eqref{def:E1234} in Appendix~\ref{appA}, $G-G_0$ decomposes into the sum of four matrices, where it is most subtle to bound the spectral norm of the fourth matrix:
\[
E_4:= M_0^{-1/2}(ZZ'-\mathbb{E}[ZZ'])M_0^{-1/2}. 
\]

Define $X_i=(M_0^{-1/2}z_i)(M_0^{-1/2}z_i)' - \mathbb E  [(M_0^{-1/2}z_i \big) \big(M_0^{-1/2}z_i)']$. It is seen that $E_4=\sum_{i=1}^n X_i$, which is a sum of $n$ independent matrices. We apply the matrix Bernstein inequality \cite{tropp2012user} (Theorem~\ref{thm:matrixBernstein}) to obtain that
if there exist $b>0$ and $\sigma^2>0$ such that $\Vert X_i \Vert\leq b$ almost surely for all $i$ and $ \Vert \sum_{i=1}^n \mathbb{E} X_i^2\Vert\leq \sigma^2$, then for every $t>0$, 
\[
\mathbb{P} \Bigl( \big \Vert\sum_{i=1}^n X_i \big\Vert \geq t \Bigr)\leq 2p \exp\Big( -\frac{t^2/2}{\sigma^2 + bt/3}\Big).
\] 

Determination of $b$ and $\sigma^2$ requires upper bounds for $\|X_i\|$ and $\|\mathbb E  X_i^2 \|$. Since each $X_i$ is equal to a rank-1 matrix minus its expectation, it reduces to deriving large-deviation bounds for $\|M_0^{-1/2}z_i\|^2$.  Note that each $z_i$ can be equivalently represented by $z_i = N_i^{-1} \sum_{m=1}^N (T_{im} - \mathbb E T_{im})$, where $\{T_{im}\}_{m=1}^{N_i}$ are i.i.d. Multinomial$(1, d_i^0)$. 
It yields that $\|M_0^{-1/2}z_i\|^2={\cal I}_1+{\cal I}_2$, where ${\cal I}_2$ is a term that can be controlled using standard large-deviation inequalities, and:
\[
{\cal I}_1:  = N_i^{-2} \sum_{1\leq m_1\neq m_2\leq N_i}  (T_{im_1} - \mathbb E T_{im_1}) M_0^{-1} (T_{im_2} - \mathbb E T_{im_2}). 
\]

The remaining question is how to bound $|{\cal I}_1|$. We notice that $\mathcal I_1$ is a U-statistic with degree~2. 
The decoupling inequality \citep{de1995decoupling, de2012decoupling} is a useful tool for studying U-statistics.
\begin{Theorem}[A special decoupling inequality \citep{de2012decoupling}]\label{thm:main_decoupling2}
Let $\{X_{m}\}_{m=1}^N$ be a sequence of i.i.d. random vectors in $\mathbb{R}^{d}$, and let $\{\widetilde X_m\}_{m=1}^N$ be an independent copy of $\{X_{m}\}_{m=1}^N$. Suppose that $h: \mathbb R^{2d}\to \mathbb R$ is a measurable function. Then, there exists a constant $C_4>0$ independent of $n, m, d$ such that for all $t>0$:
\begin{align*}
\mathbb P \Big( \Big| \sum_{m\neq m_1} h (X_m, X_{m_1}) \Big| \geq t  \Big) \leq C_4 \mathbb P \Big( C_4\Big| \sum_{m\neq m_1} h (X_m, \widetilde X_{m_1}) \Big| \geq t  \Big).  
\end{align*}
\end{Theorem}

Let $\{\widetilde T_{im}\}_{m=1}^{N_i}$ be an independent copy of $\{T_{im}\}_{m=1}^{N_i}$.
By Theorem~\ref{thm:main_decoupling2}, the large-deviation bound of $\mathcal I_1$ can be inferred from the large-deviation bound of:
\[
\widetilde{\mathcal I}_1: = N_i^{-2} \sum_{1\leq m_1\neq m_2\leq N_i} (T_{im_1} - \mathbb E T_{im_1})' M_0^{-1} (\widetilde T_{im_2} - \mathbb E \widetilde T_{im_2})\,.
\]

Using $h(T_{im_1},\widetilde T_{im_2})$ to denote the summand in the above sum, we have a decomposition: $\widetilde{\cal I}_1=N_i^{-2}\sum_{m_1,m_2}h(T_{im_1},\widetilde T_{im_2})- N_i^{-2}\sum_{m}h(T_{im},\widetilde T_{im})$. The second term is a sum of independent variables and can be controlled by standard large-deviation inequalities. Hence, the analysis of $\widetilde{\cal I}_1$ reduces to the analysis of $\widetilde{\mathcal I}_1^*:=N_i^{-2}\sum_{m_1,m_2}h(T_{im_1}, \widetilde T_{im_2})$. We re-write $\widetilde{\mathcal I}_1^*$ as:
\[
\widetilde{\mathcal I}_1^* = N_i^{-2}y'\tilde{y}, \quad\mbox{with}\;\; y:=\sum_{m=1}^{N_i} M_0^{-1/2}(T_{im} - \mathbb E T_{im}), \quad \tilde{y}:=\sum_{m=1}^{N_i} M_0^{-1/2}(\widetilde{T}_{im} - \mathbb E \widetilde{T}_{im}). 
\]

Since $\tilde{y}$ is independent of $y$, we apply large-deviation inequalities {twice}. First, conditional on $\tilde{y}$, $\widetilde{\mathcal I}_1^*$ is a sum of $N_i$ independent variables (randomness comes from $T_{im}$'s). We apply the Bernstein inequality to get a large-deviation bound for $\widetilde{\cal I}_1^*$, which depends on a quantity $\sigma^2(\tilde{y})$. Next, since $\sigma^2(\tilde{y})$ can also be written as a sum of $N_i$ independent variables (randomness comes from $\widetilde{T}_{im}$'s), we apply the Bernstein inequality again to obtain a large-deviation bound for $\sigma^2(\tilde{y})$. 
Combining two steps gives the large-deviation bound for $\widetilde{\mathcal I}_1^*$.

\begin{Remark} \label{rmk:decoupling}
The decoupling inequality is employed multiple times to study other U-statistics-type quantities arising in our proof. For example,  recall that $(G-G_0)$ decomposes into the sum of four matrices, and we have only discussed how to bound $\|E_4\|$. In the analysis of $\|E_2\|$ and $\|E_3\|$, we need to bound other quadratic terms involving a sum over $(i,m)$, with $1\leq i\leq n$ and $1\leq m\leq N_i$. In that case, we need a more general decoupling inequality. We refer readers to  Theorem~\ref{thm:decoupling2} in Appendix~\ref{appA} for more details.
\end{Remark}

\begin{Remark} \label{rmk:comparison}
The analysis in \cite{ke2022using} uses an $\epsilon$-net argument \citep{Vershynin} and the martingale Bernstein inequality~\citep{freedman1975tail} to study $\|E_4\|$. In our analysis, we use the matrix Bernstein inequality \citep{tropp2012user}, instead of the $\epsilon$-net argument. The matrix Bernstein inequality enables us to tackle each quadratic term related to each $i$ separately instead of handling complicated quadratic terms involving summation over $i$ and $m$ simultaneously. 
Additionally, we adopt the decoupling inequality for U-statistics \citep{de1995decoupling, de2012decoupling}, instead of the martingale Bernstein inequality, to study all the quadratic terms arising in our analysis. The decoupling inequality converts the tail anaysis of quadratic terms into tail analysis of (conditionally) independent sums. It provides sharper bounds when the variables have heavy tails (which is the case for the word counts in a topic model, especially when documents are short). 
\end{Remark}

\subsection{Proof sketch of Theorem~\ref{thm:row_eigenv}}

We combine the non-stochastic perturbation result in Lemma~\ref{lem:perturb} and the large-deviation bounds in Lemma~\ref{lem:large-dev-G} to prove Theorem~\ref{thm:row_eigenv}. By Lemma~\ref{lem:eigen}, $
|\lambda_K|\geq C^{-1} n\beta_n$. 
It follows from Weyl's inequality, the first claim in Lemma~\ref{lem:large-dev-G}, and the assumption of $p \log^2 (n) \leq Nn\beta^2_n$ that with probability $1-n^{-C_1}$:
\[
|\hat{\lambda}_K|\geq |\lambda_k|\cdot \bigl[1 - O\bigl([\log(n)]^{-1/2}\bigr)\bigr] \geq C^{-1}n\beta_n.
\] 

In addition, it  can be shown (see Lemma~\ref{lem:eigen}) that $\|e_j'\Xi\|\leq Ch_j^{1/2}$. Combining this with the two claims in Lemma~\ref{lem:large-dev-G} gives that with probability $1-n^{-C_1}$:
\[
\|e_j'\Xi\|\|G-G_0\|+\|e_j'(G-G_0)\Xi\|\leq C\sqrt{\frac{h_jnp\log(n)}{N}}:=\gamma_j. 
\] 

We hope to apply Lemma~\ref{lem:perturb}. This requires obtaining a bound for $\|\Gamma^{-1}(G-G_0)\Gamma\|_{1\to\infty}$. 
Since 
$\Gamma\propto H^{1/2}$, it suffices to study $\|H^{-1/2}(G-G_0)H^{1/2}\|_{1\to\infty}$. Similar to the analysis of $\|e_j'(G-G_0)\Xi\|$, we can show (see the proofs of Lemmas~\ref{lem:tech2} and \ref{lem:tech3}, such as \eqref{2023101001}) that $\|e_j'(G-G_0)H^{1/2}\|_1\leq CN^{-1/2}[h_jnp\log(n)]^{1/2}\leq C\sqrt{h_j/\log(n)}\cdot n\beta_n$, where the last inequality is because of $p\log^2(n)\leq Nn$. We immediately have:
\[
\|H^{-1/2}(G-G_0)H^{1/2}\|_{1\to\infty}=\max_j \bigl\{h_j^{-1/2}\|e_j'(G-G_0)H^{-1/2}\|_1\bigr\}\leq \frac{Cn\beta_n}{\sqrt{\log(n)}}\leq \frac{|\hat{\lambda}_K|}{2}. 
\]

We then apply Lemma~\ref{lem:perturb} to get $\|e_j'(\hat\Xi-\Xi O')\|\leq C|\hat{\lambda}_K|^{-1}\gamma_j\leq C(n\beta_n)^{-1}\gamma_j$. The claim of Theorem~\ref{thm:row_eigenv} follows immediately by plugging in the value of $\gamma_j$ as given above.

\section{Summary and Discussion} \label{sec:discuss}

The topic model imposes a ``low-rank plus noise'' structure on the data matrix. However, the noise is not simply additive; rather, it consists of centered multinomial random vectors. The eigenvector analysis in a topic model is more complex than standard eigenvector analysis for random matrices. Firstly, the entries of the data matrix are weakly dependent, making techniques such as leave-one-out inapplicable. Secondly, due to the significant word frequency heterogeneity in natural languages, entry-wise eigenvector analysis becomes much more nuanced, as different entries of the same eigenvector have significantly different bounds. Additionally, the data exhibit Bernstein-type tails, precluding the use of random matrix theory tools that assume sub-exponential entries. While we build on the analysis in \cite{ke2022using}, we address these challenges with new proof ideas. Our results provide the most precise eigenvector analysis and rate-optimality results for topic modeling, to the best of our knowledge.

A related but more ambitious goal is obtaining higher-order expansions of the empirical singular vectors. Since the random matrix under study in the topic model is the Wishart type, we can possibly borrow techniques in \citep{bloemendal2016principal}
to study the joint distribution of empirical singular values and singular vectors. In this paper, we assume the number of topics, $K$, is finite, but our analysis can be easily extended to the scenario of a growing $K$ (e.g., $K=O(\log(n))$). 
We assume $\min\{p, N\} \geq \log^3 (n)$. When $p<\log^3(n)$, it becomes a low-dimensional eigenvector analysis problem, which is easy to tackle. When $N<\log^3(n)$, it is the {\it extremely short documents} case (i.e., each document has only a finite length, say, fewer than $20$, as in documents such as Tweets). We leave it to future work.

\appendix


\section{Preliminary Lemmas and Theorems}\label{appA}
In this section, we collect the preliminaries lemmas and theorems that will be used in the entry-wise eigenvector analysis. Under Assumption~\ref{asm:para}, $N_i \asymp \bar N \asymp N$. Therefore, throughout this section and subsequent sections, we always assume $\bar{N}=N$ without loss of generality.

The first lemma describes the estimates of the entries in $M_0$ and reveals its relation to the underlying frequency parameters, and further provides the large-deviation bound for the normalization  matrix $M$. 
\begin{Lemma}[Lemmas D.1 \& E.1 in \citep{ke2022using}]\label{lem:MM0}
Recall the definitions $M= {\rm diag}(n^{-1} \sum_{i=1}^n N d_i/N_i) $, $M_0= {\rm diag}(n^{-1}\sum_{i=1}^n N d_i^0/N_i) $, and $h_j = \sum_{k=1}^K A_k(j)$ for $1\leq j \leq p$. Suppose the conditions in Theorem~\ref{thm:row_eigenv} hold. Then:
\begin{align*}
M_0(j,j) \asymp h_j;\qquad \text{ and} \quad   | M(j,j) - M_0(j,j) | \leq C\sqrt{\frac{h_j \log (n)}{Nn}},
\end{align*}
for some constant $C>0$, with probability $1- o(n^{-3})$, simultaneously for all $1\leq j \leq p$. Furthermore, with probability $1- o(n^{-3})$, 
\begin{align} \label{est:MM01/2}
\Big \Vert M^{-1/2} M_0^{1/2} - I_p \Big \Vert  \leq C\sqrt{\frac{p\log (n)}{Nn}}.
\end{align}
\end{Lemma}

\begin{Remark}
In this lemma and other subsequent lemmas, ``with probability $1-o(n^{-3})$'' can always be replaced by ``with probability $1-n^{-C_1}$'', for an arbitrary constant $C_1>0$. The small-probability events in these lemmas come from the Bernstein inequality or the matrix Bernstein inequality. These inequalities concern small-probability events associated with an arbitrary probability $\delta\in (0,1)$, and the high-probability bounds depend on $\log(1/\delta)$. When $\delta=n^{-C_1}$, $\log(1/\delta)=C_1\log(n)$. Therefore, changing $C_1$ only changes the high-probability bound by a constant. Without loss of generality, we take $C_1=4$ for convenience. 
\end{Remark}

The proof of the first statement is quite similar to the proof detailed in the supplementary materials of \citep{ke2022using}. The only difference is the existence of the additional factor $N/N_i$. Thanks to the condition that $N_i$'s are at the same order, it is not hard to see that $M_0(j,j) \asymp n^{-1}\sum_{i=1}^n d_i^0(j) $,where the RHS is exactly the definition of $M_0$ in \citep{ke2022using}. Thus, the proof follows simply under Assumption~\ref{asm:WA}. To obtain the large-deviation bound, the following representation is crucial:
\[
M(j,j) - M_0(j,j) = \frac 1n \sum_{i=1}^n \frac{N}{N_i} \big(d_i(j) - d_i^0(j) \big)= \frac{1}{n} \sum_{i=1}^n \frac{N}{N_i^2}\sum_{m=1}^{N_i} T_{im}(j)-  d_i^0(j) , 
\] 
where $\{T_{im}\}_{m=1}^n$ are i.i.d. Multinomial $(1, d_i^0)$ with $d_i^0 = Aw_i$. The RHS is a sum of independent random variables, thus allowing the application of Bernstein inequality. The inequality (\ref{est:MM01/2}) is not provided in the supplementary materials of \citep{ke2022using}, but it follows easily from the first statement. We prove (\ref{est:MM01/2})  in detail below. 

By definition, it suffices to claim that:
\[
\bigg|\frac{\sqrt {M_0(j,j)}}{ \sqrt{M(j,j)} }  - 1 \bigg| \leq C\sqrt{\frac{p\log (n)}{Nn}}
\]
simultaneously for all $1\leq j\leq p$. To this end, we derive:
\[
\bigg|\frac{\sqrt {M_0(j,j)}}{ \sqrt{M(j,j)} }  - 1 \bigg| \leq \frac{ \big| M_0(j,j) - M(j,j) \big| }{ \sqrt{M(j,j)} (\sqrt {M_0(j,j)} + \sqrt{M(j,j)}) } 
\]

Using the large-deviation bound $ | M(j,j) - M_0(j,j) | \leq C\sqrt{h_j \log (n)/({Nn})} = o (h_j)$ and also the estimate $M_0(j,j)\asymp h_j$, we bound the denominator by:
\[
\sqrt{M(j,j)} \Big(\sqrt {M_0(j,j)} + \sqrt{M(j,j)} \Big) \geq  C \sqrt{h_j - o(h_j) } \Big( \sqrt{h_j} +  \sqrt{h_j - o(h_j) } \Big) \geq Ch_j 
\]
with probability $1- o(n^{-3})$, simultaneously for all $1\leq j \leq p$. Consequently:
\[
\bigg|\frac{\sqrt {M_0(j,j)}}{ \sqrt{M(j,j)} }  - 1 \bigg| \leq C\sqrt{\frac{\log (n)}{Nn h_j}}\leq C\sqrt{\frac{p\log (n)}{Nn}}, 
\]
where the last step is due to $h_j \geq h_{\min} \geq C /p$. This completes the proof of (\ref{est:MM01/2}).

The next Lemma presents the eigen-properties of the population data matrix. 
\begin{Lemma}[Lemmas F.2, F.3, and D.3 in \citep{ke2022using}] \label{lem:eigen}
Suppose the conditions in Theorem~\ref{thm:row_eigenv} hold.  Let $G_0 $ be as in (\ref{def:GG0}). Denote by $\lambda_1\geq \lambda_1\geq \ldots \geq \lambda_K$ the non-zero eigenvalues of $G_0$. There exists a constant $C>1$ such that:
\[
Cn \beta_n \leq \lambda_k\leq C n, \;\; \text{ for $2\leq k\leq K$},\qquad \text{and} \quad \lambda_1 \geq C^{-1} n + \max_{2\leq k \leq K}\lambda_K\,. 
\]
Furthermore, let $\xi_1,\xi_2, \ldots, \xi_K $ be the associated eigenvectors of $G_0$. Then:
\begin{align*}
C^{-1} \sqrt{h_j} \leq \xi_1(j)\leq C \sqrt{h_j}\, , \qquad \Vert e_j'\Xi\Vert \leq C\sqrt{h_j} \, . 
\end{align*}
\end{Lemma}

The above lemma can be proved in the same manner as those in the supplement materials of  \citep{ke2022using}. Given our more general condition on $\Sigma_A$, which allows its smallest eigenvalue to converge to $0$ as $n\to \infty$, the results on the eigenvalues are slightly different. In out setting, only the largest eigenvalue is of order $n$ and it is well-separated from the others as  the first eigenvector of $n^{-1}G_0$ has multiplicity one, which can be claimed by using Perron's theorem and  the last inequality in Assumption~\ref{asm:WA}. For the other eigenvalues, they might be at the order of $\beta_n$ in Assumption~\ref{asm:WA}.  The details are very similar to those  in the supplement materials of  \citep{ke2022using} by adapting our relaxed condition on $\Sigma_A$, so we avoid redundant derivations here. 

Throughout the analysis, we need matrix Bernstein inequality and decoupling inequality for U-statistics. For readers' convenience, we provide the theorems below. 
\begin{Theorem}\label{thm:matrixBernstein}
Let $X_1, \cdots, X_N$ be independent, mean zero, $n\times n $ symmetric random matrices, such that $\Vert X_i \Vert\leq b$ almost surely for all $i$ and $\Vert \sum_{i=1}^N \mathbb{E} X_i^2\Vert \leq \sigma^2$. Then, for every $t\geq 0$, we have:
\begin{align*}
\mathbb{P} \left( \Big \Vert\sum_{i=1}^N X_i \Big\Vert \geq t \right)\leq 2n \exp\Big( -\frac{t^2/2}{\sigma^2 + bt/3}\Big).
\end{align*} 
\end{Theorem}

The following two theorems are special cases of Theorem 3.4.1 in \citep{de2012decoupling}, which implies that using decoupling inequality simplifies the analysis of  U-statistics to the study of sums of (conditionally) independent random variables.

\begin{Theorem}\label{thm:decoupling}
Let $\{X_i\}_{i=1}^n$ be a sequence of i.i.d. random vectors in $\mathbb{R}^{d}$, and let $\{\widetilde X_i\}_{i=1}^n$ be an independent copy of $\{X_i\}_{i=1}^n$. Then, there exists a constant $\widetilde C>0$ independent of $n, d$ such that:
\begin{align*}
\mathbb P ( | \sum_{i\neq j} X_i' X_j  | \geq t  ) \leq \widetilde C\,  \mathbb P (\, \widetilde C\,  | \sum_{i\neq j} X_i' \widetilde X_j  | \geq t )
\end{align*}
\end{Theorem}

\begin{Theorem}\label{thm:decoupling2}
Let $\{X^{(i)}_{m}\}_{i,m}$, for $1\leq i \leq n $ and $1\leq m \leq N$, be a sequence of i.i.d. random vectors in $\mathbb{R}^{d}$, and let $\{\widetilde X^{(i)}_m\}_{i,m}$ be an independent copy of $\{X^{(i)}_{m}\}_{i,m}$. Suppose that $h: \mathbb R^{2d}\to \mathbb R$ is a measurable function. Then, there exists a constant $\overline C>0$ independent of $n, m, d$ such that:
\begin{align*}
\mathbb P \Big( \Big|\sum_{i} \sum_{m\neq m_1} h (X_m^{(i)}, X_{m_1}^{(i)}) \Big| \geq t  \Big) \leq \overline C\, \mathbb P \Big( \, \overline C\,  \Big|\sum_{i} \sum_{m\neq m_1} h (X_m^{(i)}, \widetilde X_{m_1}^{(i)}) \Big| \geq t  \Big) 
\end{align*}
\end{Theorem}

The key difference between the above theorems is attributed to the index set used across the sum. In Theorem~\ref{thm:decoupling}, the random variables are indexed by $i$ and all pairs of $(X_i, X_j)$ are included; in contrast, Theorem ~\ref{thm:decoupling2} uses both $i$ and $m$ and consider only the pairs that share the identical index $i$.  However, both are viewed as special cases of  Theorem 3.4.1 with degree $2$ in \citep{de2012decoupling}, which discussed a broader sequence of functions $\{h_{ij}(\cdot, \cdot)\}_{i,j}$, where each  $h_{ij}(\cdot, \cdot)$ can differ with varying $i,j$. By assigning all $h_{ij}(\cdot, \cdot)$ to the same product function,  we have Theorem~\ref{thm:decoupling};  whereas Theorem~\ref{thm:decoupling2} follows from specifying:
\[
h_{(im)(jm_1)} (\cdot, \cdot) = \left \{
\begin{array}{cl}
h(\cdot, \cdot), & {\text if }\quad  i=j;\\
0,&  {\text otherwise}. 
\end{array}
\right.
\]

\section{Proofs of Lemmas~\ref{lem:perturb} and \ref{lem:large-dev-G}}\label{appB}

\subsection{Proof of Lemma~\ref{lem:perturb}}
Using the definition of eigenvectors and eigenvalues, we have $G\hat\Xi=\hat\Xi\hat\Lambda$ and $G_0\Xi=\Xi\Lambda$. Additionally, since $G_0$ has a rank $K$, $G_0=\Xi\Lambda\Xi'$. It follows that:
\[
\hat \Xi  \hat \Lambda =  [G_0 + (G- G_0)]\hat  \Xi = \Xi\Lambda \Xi' \hat \Xi + (G-G_0) \hat\Xi = \Xi \Xi'G_0 \hat \Xi  + (G-G_0) \hat\Xi. 
\]
As a result:
\begin{equation} \label{proof-perturb-1}
e_j' \hat \Xi = e_j' \Xi \Xi'G_0 \hat \Xi \hat \Lambda^{-1} + e_j'(G- G_0) \hat \Xi \hat \Lambda^{-1}\,. 
\end{equation}
Note that $G_0\hat{\Xi}=G\hat{\Xi}+(G_0-G)\hat{\Xi}=\hat\Xi\hat{\Lambda}+(G_0-G)\hat\Xi$. We plug this equality into the first term on the RHS of \eqref{proof-perturb-1} to obtain: 
\begin{align*}
 e_j' \Xi \Xi'G_0 \hat \Xi \hat \Lambda^{-1}  &= e_j' \Xi\Xi' \hat \Xi + e_j' \Xi \Xi' (G_0- G) \hat \Xi \hat \Lambda^{-1}\notag\\
& = e_j' \Xi O' + e_j' \Xi  (\Xi ' \hat \Xi - O') + e_j' \Xi \Xi' (G_0- G) \hat \Xi \hat \Lambda^{-1},
\end{align*}
for any orthogonal matrix $O$. Combining this with \eqref{proof-perturb-1} gives:
\begin{equation}\label{proof-perturb-2}
\Vert e_j' (\hat \Xi- \Xi O')   \Vert  \leq \Vert  e_j' \Xi (\Xi' \hat \Xi - O') \Vert  + \Vert e_j' \Xi \Xi' (G_0- G) \hat \Xi  \hat \Lambda^{-1}\Vert +  \Vert e_j'(G- G_0) \hat \Xi \hat \Lambda^{-1} \Vert.
\end{equation}
Fix $O={\rm sgn} (\hat \Xi' \Xi)$. 
The sine-theta theorem \citep{sin-theta} yields:
\begin{equation} \label{proof-perturb-3}
\|\Xi'\hat\Xi - O'\|\leq |\hat{\lambda}_K|^{-2} \|G-G_0\|^2. 
\end{equation}
We use \eqref{proof-perturb-3} to bound the first two terms on the RHS of \eqref{proof-perturb-2}:
\begin{align*}
\Vert  e_j' \Xi (\Xi' \hat \Xi - O') \Vert & \leq \|e_j'\Xi\|\|\Xi'\hat\Xi - O'\| \leq \|e_j'\Xi\|\cdot |\hat{\lambda}_K|^{-2} \|G-G_0\|^2,\cr
\Vert e_j' \Xi \Xi' (G_0- G) \hat \Xi  \hat \Lambda^{-1}\Vert &\leq \|e_j'\Xi\|\cdot  |\hat{\lambda}_K|^{-1}\| \Xi' (G_0- G) \hat \Xi \|\leq \|e_j'\Xi\|\cdot  |\hat{\lambda}_K|^{-1} \| G- G_0\|. 
\end{align*}
Since $\| G- G_0\|\leq (1-c_0)|\hat{\lambda}_K|$, the RHS in the second line above dominates the RHS in the first line. We plug these upper bounds into \eqref{proof-perturb-2} to get:
\begin{align} \label{proof-perturb-4}
\Vert e_j' (\hat \Xi- \Xi O')   \Vert  & \leq |\hat{\lambda}_K|^{-1}\|e_j'\Xi\| \| G- G_0\| + \Vert e_j'(G- G_0) \hat \Xi \hat \Lambda^{-1} \Vert\cr
&\leq |\hat{\lambda}_K|^{-1}\bigl( \|e_j'\Xi\| \| G- G_0\| + \Vert e_j'(G- G_0) \hat \Xi \Vert\bigr). 
\end{align}

We notice that the second term on the RHS of \eqref{proof-perturb-4} still involves $\hat{\Xi}$, and we further bound this term. 
By the assumption of this theorem, there exists a diagonal matrix $\Gamma$ such that $\|\Gamma^{-1}(G-G_0)\Gamma\|_{1\to\infty}\leq (1-c_0)|\hat{\lambda}_K|$. It implies:
\[
\|e_j'(G-G_0)\Gamma\|_1 \leq (1-c_0)\gamma_j|\hat\lambda_K|.
\]
Additionally, for any vector $v\in\mathbb{R}^p$ and matrix $B\in\mathbb{R}^{p\times K}$, it holds that $\|v'B\|\leq \sum_j |v_j|\|e_j'B\|\leq \sum_j |v_j|\|B\|_{2\to\infty}\leq \|v\|_1\|B\|_{2\to\infty}$. 
We then bound the second term on the RHS of \eqref{proof-perturb-4} as follows:
\begin{align} \label{proof-perturb-5}
\Vert e_j'(G- G_0) \hat \Xi \Vert & \leq \|e_j'(G-G_0)\Xi O'\| + \|e_j'(G-G_0)(\hat\Xi - \Xi O')\|\cr
&\leq \|e_j'(G-G_0)\Xi \| +\|e_j'(G-G_0)\Gamma\|_1\cdot\|\Gamma^{-1}(\hat\Xi-\Xi O')\|_{2\to\infty}\cr
&\leq \|e_j'(G-G_0)\Xi \| + (1-c_0)\gamma_j|\hat{\lambda}_K|\cdot\|\Gamma^{-1}(\hat\Xi-\Xi O')\|_{2\to\infty}. 
\end{align}
Plugging \eqref{proof-perturb-5} into \eqref{proof-perturb-4} gives:
\begin{align}  \label{proof-perturb-6}
\Vert e_j' (\hat \Xi- \Xi O')   \Vert &\leq |\hat{\lambda}_K|^{-1} \bigl( \|e_j'\Xi\| \| G- G_0\| + \|e_j'(G-G_0)\Xi \| \bigr) \cr
&\quad + (1-c_0)\gamma_j\cdot \Vert \Gamma^{-1}(\hat\Xi-\Xi O')\Vert_{2\to\infty}\cr
&\leq |\hat{\lambda}_K|^{-1}\gamma_j + (1-c_0)\gamma_j\cdot \Vert \Gamma^{-1}(\hat\Xi-\Xi O')\Vert_{2\to\infty},
\end{align}
where in the last line we have used the assumption that $\gamma_j$ is an upper bound for $\|e_j'\Xi\| \| G- G_0\| + \|e_j'(G-G_0)\Xi \|$. 
Note that $\Vert \Gamma^{-1}(\hat\Xi-\Xi O')\Vert_{2\to\infty}= \max_{1\leq j\leq p}\bigl\{\gamma_j^{-1}\Vert e_j' (\hat \Xi- \Xi O')   \Vert\bigr\}$. We multiply both LSH and RSH of \eqref{proof-perturb-6} by $\gamma_j^{-1}$  and take the maximum over $j$. It gives:
\begin{equation} \label{proof-perturb-7}
\Vert \Gamma^{-1}(\hat\Xi-\Xi O')\Vert_{2\to\infty} \leq |\hat{\lambda}_K|^{-1} +  (1-c_0)\Vert \Gamma^{-1}(\hat\Xi-\Xi O')\Vert_{2\to\infty}, 
\end{equation}
or equivalently, $\|\Gamma^{-1}(\hat\Xi-\Xi O')\|_{2\to\infty} \leq c^{-1}_0|\hat\lambda_K|^{-1}$. We further plug this inequality into~\eqref{proof-perturb-6} to obtain:
\begin{equation} \label{proof-perturb-8}
\|e_j' (\hat \Xi- \Xi O')\|\leq |\lambda_K|^{-1}\gamma_j + (1-c_0)\cdot c_0^{-1}|\lambda_K|^{-1}\gamma_j \leq c_0^{-1}|\lambda_K|^{-1}\gamma_j. 
\end{equation}
This proves the claim.

\subsection{Proof of Lemma~\ref{lem:large-dev-G}}
The first claim is the same as the one in Lemma~\ref{lem:G-G0} and will be proved there. 

The second claim follows by simply collecting arguments in the proof of Lemma~\ref{lem:G-G0}, as shown below: 
By \eqref{eq:G-G_0}, $G-G_0=E_1+E_2+E_3+E_4$. It follows that:
\begin{equation} \label{proof-largeDevG-0}
\|e_j'(G-G_0)\Xi\|\leq \sum_{s=1}^4\|e_j'E_s\Xi\|. 
\end{equation}

We apply Lemma~\ref{lem:tech2} to get large-deviation bounds for $\|e_j'E_s\Xi\|$ with $s\in\{2,3,4\}$. This lemma concerns $\|e_j'E_s\hat{\Xi}\|$, but in its proof we have already analyzed $\|e_j'E_s\Xi\|$. In particular, $\|e_j'E_2\Xi\|$ and $\|e_j'E_3\Xi\|$ have the same bounds as in \eqref{eq:entryE23}, and the bound for $\|e_j'E_4\Xi\|$ only has the first term in \eqref{eq:tech2_entryE4}. In summary: 
\begin{equation} \label{proof-largeDevG-1}
\|e_j'E_s\Xi\| \leq  C\sqrt{\frac{h_j np\log (n)}{N}}, \qquad\mbox{for }s\in \{2,3,4\}. 
\end{equation}

It remains to bound $\|e_j'E_1\Xi\|$. We first mimic the steps of proving (\ref{eq:tech3_entryE1}) of Lemma~\ref{lem:tech2} (more specifically, the derivation of \eqref{lem-A6-quote}, except that $\hat{\Xi}$ is replaced by $\Xi$) to obtain:
\begin{align} \label{proof-largeDevG-2}
\Vert e_j E_1 \Xi \Vert &\leq Cn\Vert e_j' (M_0^{1/2} M^{-1/2}-I_p) \Xi \Vert + C\Vert e_j' G_0 (M_0^{1/2}M^{-1/2}-I_p)\Xi  \Vert   \notag\\
& \qquad + \sum_{s=2}^4  \Vert e_j' E_s (M_0^{1/2}M^{-1/2}-I_p) \Xi  \Vert.
\end{align} 
We note that:
\begin{align*}
\Vert e_j' (M_0^{1/2} M^{-1/2}-I_p) \Xi \Vert & \leq \|M_0^{1/2} M^{-1/2}-I_p\|\cdot \|e_j'\Xi\|,\cr
\Vert e_j' G_0 (M_0^{1/2}M^{-1/2}-I_p)\Xi  \Vert &= \Vert e_j' \Xi\Lambda \Xi' (M_0^{1/2}M^{-1/2}-I_p)\Xi  \Vert \cr
&\leq \|e_j'\Xi\|\cdot\|\Lambda\|\cdot \|M_0^{1/2} M^{-1/2}-I_p\|,\cr
\Vert e_j' E_s (M_0^{1/2}M^{-1/2}-I_p) \Xi\| &\leq \|e_j'E_s\|\cdot  \|M_0^{1/2} M^{-1/2}-I_p\|.  
\end{align*}
For $s\in \{2,3\}$, we have $\|e_j'E_s\|\leq C\sqrt{h_jp\log(n)/(Nn)}$. This has been derived in the proof of Lemma~\ref{lem:tech2}: when controlling $\|e_j'E_2\Xi\|$ and $\|e_j'E_3\Xi\|$ there, we first bound them by $\|e_j'E_2\|$ and $\|e_j'E_3\|$, respectively, and then study $\|e_j'E_2\|$ and $\|e_j'E_3\|$ directly). We plug these results into \eqref{proof-largeDevG-2} to obtain:
\begin{align} \label{proof-largeDevG-3}
\Vert e_j E_1 \Xi \Vert &\leq \|M_0^{1/2}M^{-1/2}-I_p\|\Biggl(n\|e_j'\Xi\| +|\lambda_1| \|e_j'\Xi\|+  C\sqrt{\frac{h_j n p\log (n)}{N}} \Biggr)\cr
&\qquad + \Vert e_j' E_4 (M_0^{1/2}M^{-1/2}-I_p) \Xi\|. 
\end{align}
For $\|e_j'E_4(M_0^{1/2}M^{-1/2}-I_p)\Xi\|$, we cannot use the same idea to bound it as for $s\in\{2,3\}$, because the bound for $\|e_j'E_4\|$ is much larger than those for $\|e_j'E_2\|$ and $\|e_j'E_4\|$. Instead, we study $\|e_j'E_4(M_0^{1/2}M^{-1/2}-I_p)\Xi\|$ directly. This part is contained in the proof of Lemma~\ref{lem:tech3}; specifically, in the proof of 
\eqref{eq:tech2_add}. There we have shown:
\begin{equation} \label{proof-largeDevG-4}
\Vert e_j' E_4 (M_0^{1/2}M^{-1/2}-I_p) \Xi\|\leq C\sqrt{h_j}\cdot \frac{p\log(n)}{N}. 
\end{equation}
We plug \eqref{proof-largeDevG-4} into \eqref{proof-largeDevG-3} and note that $\lambda_1=O(n)$ and $\|e_j'\Xi\|=O(h_j^{1/2})$ (by Lemma~\ref{lem:eigen}). We also use the assumption that $Nn\geq Nn\beta_n^2 \geq p\log^2(n)$ and the bound for $\|M_0^{1/2}M^{-1/2}-I_p\|$ in \eqref{est:MM01/2}. It follows that
\begin{align} \label{proof-largeDevG-5}
\Vert e_j E_1 \Xi \Vert& \leq \Vert M_0^{1/2}M^{-1/2}-I_p\Vert\cdot C\sqrt{h_j}\Biggl( n+
\sqrt{\frac{np\log(n)}{N}}+ \frac{p\log(n)}{N}\Biggr)\cr
&\leq \Vert M_0^{1/2}M^{-1/2}-I_p\Vert\cdot O(nh_j^{1/2})\;\; \leq\;\; C\sqrt{\frac{h_jnp\log(n)}{N}}\ . 
\end{align}

We plug \eqref{proof-largeDevG-1} and \eqref{proof-largeDevG-5} into \eqref{proof-largeDevG-0}. This proves the second claim.

\section{The complete proof of Theorem~\ref{thm:row_eigenv}} \label{sec:mainproof}

A proof sketch of Theorem~\ref{thm:row_eigenv} has been given in Section~\ref{subsec:U-statistics}. 
For the ease of writing formal proofs, we have re-arranged the claims and analyses in Lemmas~\ref{lem:perturb} and \ref{lem:large-dev-G}, so the proof structure here is slightly different from the sketch in Section~\ref{subsec:U-statistics}. For example, Lemma~\ref{lem:G-G0} combines the claims of Lemma~\ref{lem:large-dev-G} with some steps in proving Lemma~\ref{lem:perturb}; the remaining steps in the proof of Lemma~\ref{lem:perturb} are combined into the proof of the main theorem.

First, we present a key technical lemma. The proof of this lemma is quite involved and relegated to {Appendix}~\ref{proofsec:lemG-G0}.

\begin{Lemma} \label{lem:G-G0}
Under the setting of Theorem \ref{thm:row_eigenv}. Recall $G, G_0$ in (\ref{def:GG0}). With probability $1- o(n^{-3})$:
\begin{align}
&\Vert G- G_0 \Vert  \leq  C\sqrt{\frac{pn\log(n)}{N}} \ll n\beta_n;   \label{est:G-G0}
  \\ 
&
\Vert e_j'(G- G_0) \hat \Xi \Vert/n \leq C \sqrt{\frac{h_j p\log (n)}{nN}}  \bigg(1+ \Vert H^{-\frac 12}  (\hat \Xi-\Xi O') \Vert_{2\to \infty} \bigg)+ o(\beta_n) \cdot \Vert e_j' (\hat \Xi-\Xi O')\Vert \, , \label{est:jG-G0}
\end{align}
simultaneously for all $1\leq j \leq p$. 
\end{Lemma}

Next, we use Lemma~\ref{lem:G-G0} to prove Theorem~\ref{thm:row_eigenv}. 
Let $(\hat \lambda_k, \hat \xi_k)$ and $(\hat \lambda_k, \hat \xi_k)$ be the $k$-th eigen-pairs of $G$ and $G_0$, respectively. Let $\hat \Lambda = {\rm diag} (\hat \lambda_1, \hat \lambda_2, \ldots, \hat \lambda_K)$ and $\Lambda = {\rm diag} ( \lambda_1,  \lambda_2, \ldots,  \lambda_K)$. Following \eqref{proof-perturb-1} and \eqref{proof-perturb-2}, we have:
\begin{align}\label{2022041901}
\Vert e_j' (\hat \Xi- \Xi O')   \Vert  \leq \Vert  e_j' \Xi (\Xi' \hat \Xi - O') \Vert  + \Vert e_j' \Xi \Xi' (G_0- G) \hat \Xi  \hat \Lambda^{-1}\Vert +  \Vert e_j'(G- G_0) \hat \Xi \hat \Lambda^{-1} \Vert.
\end{align}

In the sequel, we bound the three terms on the RHS above  one-by-one. 

First, by sine-theta theorem:
\begin{align*}
\Vert  e_j' \Xi (\Xi' \hat \Xi - O') \Vert \leq C \Vert e_j' \Xi \Vert \frac{\Vert G-G_0\Vert^2 }{|\hat \lambda_K- \lambda_{K+1}|^2}.
\end{align*}
For $1\leq k \leq  p $, by Weyl's inequality:
\begin{align}\label{ineq:weyl_lambda}
| \hat \lambda_k - \lambda_k| \leq \Vert G - G_0\Vert  \ll n \beta_n
\end{align}
with probability $1- o(n^{-3})$, by employing (\ref{est:G-G0}) in Lemma~\ref{lem:G-G0}. In
particular,  $\lambda_1\asymp n $ and $Cn \beta_n<\lambda_k\leq Cn$ for $2\leq k\leq K$ and $\lambda_k=0$ otherwise (see Lemma~\ref{lem:eigen}). Thereby,  $|\hat \lambda_K- \lambda_{K+1}|\geq C n \beta_n$.  Further using  $\Vert e_j' \Xi\Vert \leq C\sqrt{h_j}$ (see Lemma~\ref{lem:eigen}),
with the aid of Lemma~\ref{lem:G-G0}, we obtain that with probability $1- o(n^{-3})$:
\begin{align}\label{est:Term1}
\Vert  e_j' \Xi (\Xi'\hat \Xi - O') \Vert \leq C\sqrt{h_j} \,\cdot  \frac{p\log (n)}{Nn\beta_n^2} 
\end{align}
simultaneously for all $1\leq j \leq p$.

Next, we similarly bound the second term:
\begin{align}\label{est:Term2}
\Vert e_j' \Xi \Xi' (G_0- G) \hat \Xi \hat  \Lambda^{-1}\Vert \leq \frac{C}{n\beta_n} \Vert e_j' \Xi \Vert \Vert G- G_0\Vert  \leq C \sqrt{\frac{h_j p \log (n)}{Nn \beta_n^2}}\, .
\end{align} 
Here we used the fact that $\hat \lambda_K \geq C n\beta_n$ following from (\ref{ineq:weyl_lambda}) and Lemma~\ref{lem:eigen}.

For the last term, we simply  bound:
\begin{align}\label{est:Term3}
\Vert e_j'(G- G_0) \hat \Xi \hat \Lambda^{-1} \Vert \leq C \Vert e_j'(G- G_0) \hat \Xi \Vert/(n\beta_n)\,. 
\end{align}

Combining (\ref{est:Term1}), (\ref{est:Term2}), and (\ref{est:Term3}) into (\ref{2022041901}), by (\ref{est:jG-G0}) in Lemma~\ref{lem:G-G0},  we arrive at: 
\begin{align*}
\Vert e_j' (\hat \Xi - \Xi O')   \Vert \leq C  \sqrt{\frac{h_j p\log (n)}{Nn\beta_n^2}} \bigg(1+ \Vert H^{-\frac 12}  (\hat \Xi-\Xi O') \Vert_{2\to \infty} \bigg) + o(1) \cdot \Vert e_j' ( \hat \Xi-\Xi O')\Vert  \, .
\end{align*}
Rearranging both sides above gives: 
\begin{align}\label{eq:almost}
\Vert e_j' (\hat \Xi -  \Xi O')   \Vert \leq C  \sqrt{\frac{h_j p\log (n)}{Nn\beta_n^2}} \bigg(1+ \Vert H^{-\frac 12}  (\hat \Xi-\Xi O') \Vert_{2\to \infty} \bigg)   ,
\end{align}
with probability $ 1- o(n^{-3})$, simultaneously for all $1\leq j \leq p$. 

To proceed, we multiply both sides in (\ref{eq:almost}) by $h_j^{-1/2} $ and take the maximum. It follows that:
\begin{align*}
\Vert H^{-\frac 12}  (\hat \Xi -\Xi O') \Vert_{2\to \infty}  \leq C   \sqrt{\frac{p\log (n)}{Nn\beta_n^2}} \bigg(1 +  \Vert H_0^{-\frac 12}  (\hat \Xi -\Xi O') \Vert_{2\to \infty} \bigg)\, .
\end{align*}
Note that $\sqrt{p\log( n)}/\sqrt{Nn\beta_n^2} = o(1)$ from Assumption~\ref{asm:para}. We further rearrange both sides above and get:
\[
\Vert H^{-\frac 12}  (\hat \Xi - \Xi O') \Vert_{2\to \infty}  \leq \sqrt{\frac{p\log (n)}{Nn\beta_n^2}} = o(1)\, .
\]
Plugging the above estimate into (\ref{eq:almost}), we finally conclude 
the proof of Theorem \ref{thm:row_eigenv}.

\section{Entry-Wise Eigenvector Analysis and Proof of Lemma~\ref{lem:G-G0}}\label{appD}
To finalize the proof of Theorem~\ref{thm:row_eigenv} as outlined in {Appendix}~\ref{sec:mainproof}, the remaining task is to prove Lemma~\ref{lem:G-G0}. 

Recall the definition in (\ref{def:GG0}) that:
\[
G = M^{-\frac 12} DD' M^{-\frac 12}-\frac nN I_p,\qquad   G_0=  M_0^{-\frac 12}  \Big[\sum_{i=1}^n ( 1- N_i^{-1}) d_i^0(d_i^0)' \Big]M_0^{-\frac 12} \,. 
\] 
Write  $D = D_0 + Z$, where $Z = (z_1, z_2, \ldots, z_n)$ is a mean-zero random matrix with each $Nz_i$ being centered Multinomial $(N_i, Aw_i)$. By this representation, we decompose the perturbation matrix $G - G_0$ as follows: {
\begin{align}\label{eq:G-G_0}
 G- G_0 
&= M^{-\frac 12} DD' M^{-\frac 12} - M_0^{-\frac 12} DD' M_0^{-\frac 12}
+ M_0^{-\frac 12} \big(DD' - \sum_{i=1}^n ( 1- N_i^{-1}) d_i^0(d_i^0)'- \frac nN M_0\big) M_0^{-\frac 12} 
 \notag\\
& =( M^{-\frac 12} DD' M^{-\frac 12} - M_0^{-\frac 12} DD' M_0^{-\frac 12})    + M_0^{-\frac 12} ZD_0' M_0^{-\frac 12} + M_0^{-\frac 12} D_0Z' M_0^{-\frac 12}  \notag\\
& \quad 
+ M_0^{-\frac12} (ZZ' - \mathbb{E} ZZ') M_0^{-\frac 12}\notag\\
&=E_1 + E_2 + E_3 + E_4,
\end{align}
where:
\begin{align}\label{def:E1234}
& E_1: = M^{-\frac 12} DD' M^{-\frac 12} - M_0^{-\frac 12} DD' M_0^{-\frac 12}, \notag\\
& E_2 := M_0^{-\frac 12} ZD_0' M_0^{-\frac 12}, \qquad E_3: = M_0^{-\frac 12} D_0Z' M_0^{-\frac 12}   \notag\\
& E_4: = M_0^{-\frac12} (ZZ' - \mathbb{E} ZZ') M_0^{-\frac 12}. 
\end{align}
}
Here the second step of  (\ref{eq:G-G_0}) is due to the identity:
\begin{align*}
  \mathbb{E} (ZZ')  + \sum_{i=1}^n N_i^{-1}d_i^0(d_i^0)'  - \frac nN M_0 = 0\, ,
\end{align*}
which can be obtained by:
\[
\mathbb{E} (ZZ') = \sum_{i=1}^n \mathbb Ez_i z_i' = \sum_{i=1}^n N_i^{-2}\sum_{m, s=1}^{N_i} \mathbb E( T_{im} - \mathbb E T_{im}) ( T_{is} - \mathbb E T_{is}) ',
\]
with $\{T_{im}\}_{m=1}^N$ being  i.i.d. Multinomial $(1, Aw_i)$. 

Throughout the analysis in this section, we will frequently rewrite and use:
\begin{align}\label{eq:repre_T}
z_i = \frac {1}{N_i}\sum_{m=1}^{N_i} T_{im} - \mathbb E T_{im}
\end{align}
as it introduces the sum of independent random variables. We use the notation $d_i^0 := \mathbb E d_i = \mathbb E T_{im}= Aw_i$ for simplicity. 

By (\ref{eq:G-G_0}),  in order to prove Lemma~\ref{lem:G-G0}, it suffices to study:
\[
\|E_s\| \quad \text{and} \quad  \| e_j' E_s\hat \Xi\|/n, \qquad\text{ for $ s = 1,2, 3,4$ and $1\leq j \leq p$.} 
\]
The estimates for the aforementioned quantities are provided in the following technical lemmas, whose proofs are deferred to later sections. 

\begin{Lemma}\label{lem:tech1}
Suppose the conditions in Theorem~\ref{thm:row_eigenv} hold. There exists a constant $C>0$, such that with probability $1- o(n^{-3})$:
\begin{align}
&\| E_s\| \leq C\sqrt{\frac{pn\log (n)}{N}}, 
\qquad \text{for $s=1, 2,3$}
\label{eq:tech11}\\
& \| E_4 \|  = \| M_0^{-\frac12} (ZZ' - \mathbb{E} ZZ') M_0^{-\frac 12}\| \leq C \max\Big\{ \sqrt{\frac{pn\log (n)}{N^2}}, \frac{p\log (n)}{N}\Big\}\, . \label{eq:tech2_E4}
\end{align}
\end{Lemma}

\begin{Lemma}\label{lem:tech2}
Suppose the conditions in Theorem~\ref{thm:row_eigenv} hold. There exists a constant $C>0$, such that with probability $1- o(n^{-3})$, simultaneously for all $1\leq j \leq p$:
\begin{align}
 & \| e_j' E_s\hat \Xi\|/n \leq C\sqrt{\frac{h_jp\log (n)}{Nn}}, 
\qquad \text{for $s=2,3$} \label{eq:entryE23} \\
&  \Vert e_j' E_4 \hat \Xi  \Vert /n  \leq  C\sqrt{\frac{h_j p\log (n)}{Nn}} \,  \Big(1 +  \Vert H_0^{-\frac 12} (\hat \Xi- \Xi O')\Vert_{2\to \infty} \Big) \,, \label{eq:tech2_entryE4}
\end{align}
with $O = {\rm sgn} (\hat \Xi'\Xi )$. 
\end{Lemma}

\begin{Lemma}\label{lem:tech3}
Suppose the conditions in Theorem~\ref{thm:row_eigenv} hold. There exists a constant $C>0$, such that with probability $1- o(n^{-3})$, simultaneously for all $1\leq j \leq p$:
\begin{align}
&  \Vert e_j' E_4(M_0^{1/2}M^{-1/2}-I_p) \hat \Xi  \Vert/n \leq C\sqrt{h_j} \cdot \frac{p\log (n)}{nN} \Big( 1 +   \Vert H^{-\frac 12}  (\hat \Xi - \Xi O') \Vert_{2\to \infty}\Big) , 
\label{eq:tech2_add}\\
&\Big\Vert e_j'  \big(M^{1/2}M_0^{-1/2} - I_p\big)\hat\Xi  \Big \Vert \leq C \sqrt{\frac{\log (n)}{Nn} } + o(\beta_n) \cdot \Vert e_j'(\hat \Xi - \Xi O')\Vert; \label{eq:tech3_add}
\end{align}
and furthermore:
\begin{align}
 \| e_j' E_1\hat \Xi\|/n \leq C\sqrt{\frac{h_jp\log (n)}{Nn}} \, \Big(1 +  \Vert H_0^{-\frac 12} (\hat \Xi- \Xi O')\Vert_{2\to \infty} \Big)+ o(\beta_n) \cdot \Vert e_j'(\hat \Xi - \Xi O')\Vert\,.  \label{eq:tech3_entryE1}
\end{align}
\end{Lemma}
For proving Lemmas~\ref{lem:tech1} and \ref{lem:tech2}, the difficulty lies in showing (\ref{eq:tech2_E4}) and (\ref{eq:tech2_entryE4}) as the quantity $E_4$ involves the quadratic terms of $Z$ with its dependence on $\hat \Xi$. We overcome the hurdle by decomposing $\hat \Xi = \Xi + \hat \Xi - \Xi O'$ and employing decoupling techniques (\mbox{Theorems~\ref{thm:decoupling} and \ref{thm:decoupling2}}).  Considering the expression of $E_1$, where $DD'$ is involved, the proof of~(\ref{eq:tech3_entryE1}) in Lemma~\ref{lem:tech3} significantly rely on the estimates in Lemma \ref{lem:tech2}, together with  (\ref{eq:tech2_add}) and~(\ref{eq:tech3_add}). The detailed proofs are systematically presented in subsequent sections, following the proof of Lemma~\ref{lem:G-G0}. 

\subsection{Proof of Lemma~\ref{lem:G-G0}} \label{proofsec:lemG-G0}
We employ the technical lemmas (Lemmas~\ref{lem:tech1}--\ref{lem:tech3}) to prove Lemma~\ref{lem:G-G0}. We start with (\ref{est:G-G0}).
By the representation (\ref{eq:G-G_0}), it is straightforward to obtain that:
\[
\| G - G_0\| \leq \sum_{s=1}^4 \|E_s\| \leq C\sqrt{\frac{pn\log (n)}{N}} + C \max\Big\{ \sqrt{\frac{pn\log (n)}{N^2}}, \frac{p\log (n)}{N}\Big\}
\]
for some constant $C>0$, with probability $1- o(n^{-3})$. Under Assumption~\ref{asm:para}, it follows that:
\[
 \sqrt{\frac{pn\log (n)}{N^2}} \ll \sqrt{\frac{pn\log (n)}{N}}, \qquad \frac{p\log (n)}{N} = \sqrt{\frac{pn\log (n)}{N}} \cdot \sqrt{\frac{p \log (n)}{Nn}} \ll \sqrt{\frac{pn\log (n)}{N}}
\]
and:
\[
\sqrt{\frac{pn\log (n)}{N}} = n \cdot  \sqrt{\frac{p \log (n)}{Nn}} \ll n\, .
\]
Therefore, we complete the proof of (\ref{est:G-G0}). 

Next, we show (\ref{est:jG-G0}). Similarly, using (\ref{eq:tech11}), (\ref{eq:tech2_entryE4}), and (\ref{eq:tech3_entryE1}), we have:
\begin{align*}
\Vert e_j'(G- G_0) \hat \Xi \Vert/n &  \leq \sum_{s=1}^4 \Vert e_j'E_s \hat \Xi \Vert/n \notag\\
&\leq C\sqrt{\frac{h_jp\log (n)}{Nn}} \, \Big(1 +  \Vert H_0^{-\frac 12} (\hat \Xi- \Xi O')\Vert_{2\to \infty} \Big)+ o(\beta_n) \cdot \Vert e_j'(\hat \Xi - \Xi O')\Vert\,. 
\end{align*}
This concludes the proof of Lemma~\ref{lem:G-G0}.

\subsection{Proof of Lemma~\ref{lem:tech1}}
We examine each  $\Vert E_i\Vert $ for $i=1, 2, 3,4$. We start with the easy one, $\Vert E_2 \Vert $.
Recall $D_0 = AW$. We denote by $W_k'$ the $k$-th row of W and rewrite $W= (W_1, \cdots, W_K)'$. Similarly, we use $Z_j'$, $1\leq j\leq p$ to denote $j$-th row of $Z$. Thereby, $Z = (z_1, z_2, \ldots, z_n ) = (Z_1, Z_2, \ldots, Z_p)'$.
By the definition that $E_2=  M_0^{-1/2} ZD_0' M_0^{-1/2}$, we have:
\begin{align}\label{2024022301}
\Vert E_2 \Vert&  = \Vert M_0^{-1/2} ZW' A' M_0^{-1/2}\Vert = \Big\Vert \sum_{k=1}^K M_0^{-1/2} Z W_k\cdot  A_k'M_0^{-1/2}   \Big\Vert  \notag\\
&\leq \sum_{k=1}^K \Vert M_0^{-1/2} Z W_k\Vert  \cdot \Vert A_k'M_0^{-1/2}   \Vert. 
\end{align}
We analyze each factor in the summand:
\begin{align}\label{2024022302}
\Vert M_0^{-1/2} Z W_k\Vert^2 =\sum_{j=1}^p \frac {1}{ M_0(j,j)} (Z_j'W_k)^2, \quad  \Vert A_k'M_0^{-1/2} \Vert \asymp \Vert A_k' H^{-1} A_k\Vert^{1/2} \leq C,
\end{align}
where we used the fact that $A_k(j)\leq h_j$ for $1\leq j\leq p$. Hence, what remains is to prove a high-probability bound for each $Z_j'W_k$.  By the representation (\ref{eq:repre_T}):
\begin{align*}
Z_j'W_k = \sum_{i=1}^nz_i(j)w_i(k) =  \sum_{i=1}^n\sum_{m=1}^{N_i} N_i^{-1} w_i(k) \big( T_{im}(j) - d_i^0(j) \big). 
\end{align*}
We then apply Bernstein inequality  to the RHS above.
By straightforward computations: 
\begin{align*}
{\rm var} (Z_j'W_k)  & = \sum_{i=1}^n\sum_{m=1}^{N_i} N_i^{-2} w_i(k)^2 \mathbb E \big( T_{im}(j) - d_i^0(j)\big)^2  \notag\\
& \leq \sum_{i=1}^n  N_i^{-1}w_i(k)^2 d_i^0(j) \leq \frac{h_j n}{N},
\end{align*}
and the individual bound for each summand is $C/N$. Then, one can conclude from Bernstein inequality that  with probability $1- o(n^{-3-c_0})$:
 \begin{align}\label{eq:ZjWk}
|Z_j'W_k|\leq C\sqrt{nh_j \log (n) /N} + \log(n)/N\, . 
\end{align}
As a result, considering all $1\leq j \leq p$,  under $pn^{-c_0}\leq C$ from Assumption~\ref{asm:para}, we have:
\begin{align}\label{eq:M0ZW}
\Vert M_0^{-\frac 12} Z W_k\Vert^2 \leq C\sum_{j=1}^p h_j^{-1} \cdot \Big( \frac{nh_j \log (n)}{N}  + \frac{\log(n)^2}{N^2}\Big) \leq C \frac{np\log (n)}{N}
\end{align}
with probability $1- o(n^{-3})$.  Here, in the first step, we used $M_0(j,j)\asymp h_j $; the last step is due to the conditions $h_j \geq h_{\min} \geq C/p$ and $p\log (n) \ll Nn$. Plugging (\ref{eq:M0ZW}) and (\ref{2024022302}) into (\ref{2024022301}) gives:
\begin{align}\label{bd:E_2}
\Vert E_2 \Vert\leq  C\sqrt{\frac{np\log (n)}{N}}.
\end{align}
Furthermore, by definition, $E_3 = E_2'$ and $\Vert E_3\Vert  = \Vert E_2\Vert$. Therefore, we directly conclude the upper bound for $\Vert E_3\Vert$.

Next, we study $E_4$ and prove (\ref{eq:tech2_E4}). Notice that $M_0(j,j) \asymp h_j$ for all $1\leq j \leq p$. It suffices to prove:
\begin{align}\label{bd:new2}
 \Vert H^{-\frac 12} (ZZ'- \mathbb{E} ZZ') H^{-\frac 12} \Vert
  \leq C \max\Big\{ \sqrt{\frac{pn\log (n)}{N^2}}, \frac{p\log (n)}{N}\Big\}. 
\end{align}
We prove \eqref{bd:new2} by employing Matrix Bernstein inequality (i.e., Theorem~\ref{thm:matrixBernstein}) and decoupling techniques (i.e., Theorem~\ref{thm:decoupling}). 
First, write:
\begin{align*}
H^{-\frac 12} (ZZ'- \mathbb{E} ZZ') H^{-\frac 12} &= \sum_{i=1}^n ( H^{-\frac 12}z_i )(H^{-\frac 12}z_i )' - \mathbb{E} ( H^{-\frac 12}z_i )(H^{-\frac 12}z_i )'\notag\\
& =: n \cdot  \sum_{i=1}^n \frac 1n \big(\tilde{z}_i \tilde{z}_i'  - \mathbb{E} \tilde{z}_i \tilde{z}_i' \big) \notag\\
&=: n \cdot  \sum_{i=1}^n X_i
\end{align*}
In order to get sharp bound, we employ the truncation idea by introducing: 
\begin{align*}
\widetilde{X}_i := \frac 1n \big(\tilde{z}_i \tilde{z}_i' \mathbf{1}_{ \mathcal{E}_i} - \mathbb{E}  \tilde{z}_i \tilde{z}_i' \mathbf{1}_{\mathcal{E}_i} \big), \qquad \text{where } \quad \mathcal{E}_i:= \{\tilde{z}_i' \tilde{z}_i\leq C  p /N \},
\end{align*}
for some sufficiently large $C>0$ that depends on $C_0$ (see Assumption~\ref{asm:para}) and $ \mathbf{1}_{ \mathcal{E}_i}$ represents the indicator function.
We then have:
 \begin{align} \label{22040401}
 n\sum_{i=1}^n X_i= n\sum_{i=1}^n \widetilde{X}_i - \sum_{i=1}^n \mathbb{E} (\tilde{z}_i \tilde{z}_i' {\mathbf 1}_{\mathcal{E}_i^c} )
 \end{align}
 under the event $\bigcap_{i=1}^n \mathcal{E}_i$. We will prove the large-deviation bound of $H^{-\frac 12} (ZZ'- \mathbb{E} ZZ') H^{-\frac 12} $ in the following steps.
 \begin{itemize}
 \item [(a)]  We claim that:
 \begin{align*}
  \mathbb{P} (\bigcap_{i=1}^n \mathcal{E}_i) \leq 1- \sum_{i=1}^n \mathbb{P} (\mathcal{E}_i^c) =  1- o(n^{-(2C_0+3)})\,. 
  \end{align*}
  \item [(b)]  We claim that  under the event $\bigcap_{i=1}^n \mathcal{E}_i$:
  \[
  \Big\Vert  n\sum_{i=1}^n X_i - n\sum_{i=1}^n \widetilde{X}_i \Big\Vert = o(n^{-(C_0+1)})\, .
  \]
  \item [(c)] We aim to derive a high probability bound of $n\sum_{i=1}^n \widetilde{X}_i$ by Matrix Bernstein inequality (i.e., Theorem~\ref{thm:matrixBernstein}).  We show that with probability $1- o(n^{-3})$, for some large $C>0$:
\begin{align*}
 \Big\Vert \sum_{i=1}^n \widetilde{X}_i  \Big\Vert\leq C \max\Big\{ \sqrt{\frac{p\log (n)}{nN^2}}, \frac{p\log (n)}{nN}\Big\} \,. 
\end{align*}
 \end{itemize}
If (a)--(c) are claimed, with the condition that $N<C n^{-C_0}$ from Assumption~\ref{asm:para},  it is straightforward to conclude that:
 \begin{align*}
 \Vert H^{-\frac 12} (ZZ'- \mathbb{E} ZZ') H^{-\frac 12} \Vert  & =n \Big\|  \sum_{i=1}^n \widetilde{X}_i \Big\| + o(n^{-C_0}) \notag\\
 & \leq C \max\Big\{ \sqrt{\frac{pn\log (n)}{N^2}}, \frac{p\log (n)}{N}\Big\}\, , 
 \end{align*}
 with probability $1- o(n^{-3})$. This gives \eqref{eq:tech2_E4}, except that we still need to verify (a)--(c). 
 
 In the sequel, we prove (a), (b) and (c) separately. To prove (a), it suffices to show that $ \mathbb{P} (\mathcal{E}_i^c)=  o(n^{-(2C_0+4)})$ for all $1\leq i \leq n $.
By definition, for any fixed $i$, $N_iz_i$ is centered multinomial with $N_i$ trials. Therefore, we can represent: 
\begin{align}\label{eq:repre-z}
z_i =  \frac{1}{N_i}\sum_{m=1}^{N_i} (T_{im}- \mathbb{E} T_{im}), \quad  \text{ where $T_{im}$'s are i.i.d. multinomial$(1, d_i^0)$ for fixed $i$}, 
\end{align}
Then it can be computed that:
 \begin{align}\label{bd:expec}
 \mathbb{E} ( \tilde{z}_i' \tilde{z}_i)= \mathbb{E} z_i' H^{-1} z_i &=\frac{1}{N_i^2} \sum_{m=1}^{N_i} \mathbb{E} (T_{im}- \mathbb{E} T_{im})' H^{-1} (T_{im}- \mathbb{E} T_{im}) \notag\\
 &= \frac{1}{N_i^2} \sum_{m=1}^{N_i}\sum_{t=1}^p \mathbb{E}(T_{im}(t)- d_i^0(t))^2 h_t^{-1}\notag\\
 &=  \frac{1}{N_i^2} \sum_{m=1}^{N_i} \sum_{t=1}^p d_i^0(t)\big(1 - d_i^0(t)\big)h_t^{-1}\leq \frac{p}{N_i}.
 \end{align}
We write:
 \begin{align} \label{eq23032305}
\tilde{z}_i' \tilde{z}_i -  \mathbb{E} ( \tilde{z}_i' \tilde{z}_i) =z_i' H^{-1} z_i -  \mathbb{E} z_i' H^{-1} z_i = \mathcal{I}_1 + \mathcal{I}_2, 
 \end{align}
 where:
 \begin{align*}
& \mathcal{I}_1 :=\frac{1}{N_i^2} \sum_{m_1\neq m_2}^{N_i} (T_{im_1}- \mathbb{E} T_{im_1})' H^{-1} (T_{im_2}- \mathbb{E} T_{im_2}), \notag\\
& 
 \mathcal{I}_2:=  \frac{1}{N_i^2} \sum_{m=1}^{N_i} (T_{im}- \mathbb{E} T_{im})' H^{-1} (T_{im}- \mathbb{E} T_{im})  - \mathbb{E} (T_{im}- \mathbb{E} T_{im})' H^{-1} (T_{im}- \mathbb{E} T_{im}). 
 \end{align*}
 
First, we study ${\cal I}_1$. 
Let $\{\widetilde{T}_{im}\}_{m=1}^N $ be an independent copy of $\{{T}_{im}\}_{m=1}^N $ and:
 \begin{align*}
 \widetilde{\mathcal{I}}_1:= \frac{1}{N_i^2} \sum_{m_1\neq m_2}^{N_i} (T_{im_1}- \mathbb{E} T_{im_1})' H^{-1} (\widetilde{T}_{im_2}- \mathbb{E} \widetilde{T}_{im_2}).
 \end{align*}
We apply Theorem~\ref{thm:decoupling} to $\mathcal{I}_1$ and get:
 \begin{align}\label{eq:decoup}
 \mathbb{P}(|\mathcal{I}_1|> t) \leq C \mathbb{P}(\widetilde{\mathcal{I}}_1> C^{-1}t). 
 \end{align}
It suffices to obtain the large-deviation of $\widetilde{\mathcal{I}}_1$ instead. Rewrite:
 \begin{align}\label{form tI_1}
 \widetilde{\mathcal{I}}_1  &= \frac{1}{N_i} \sum_{m_1}^{N_i} (\widetilde{T}_{im_1}- \mathbb{E} \widetilde{T}_{im_1})' H^{-1/2} \Big( \frac {1}{N_i}  \sum_{m=1}^{N_i} H^{-1/2}(T_{im} - \mathbb{E} T_{im}) \Big)\notag\\
 & \quad - \frac{1}{N_i^2} \sum_{m=1}^{N_i}  (T_{im}- \mathbb{E} T_{im})' H^{-1} (\widetilde{T}_{im}- \mathbb{E} \widetilde{T}_{im}) \notag\\
 &=: \mathcal{T}_1 + \mathcal{T}_2. 
 \end{align}
We derive the high-probability bound for $\mathcal{T}_1$ first. For simplicity, write:
\[
a= H^{-1/2}\Big( \frac{1}{N_i} \sum_{m=1}^{N_i} ({T}_{im}- \mathbb{E} {T}_{im}) \Big) .
\]
Then, $\mathcal{T}_1 =N_i^{-1} \sum_{m=1}^{N_i} (\widetilde{T}_{im}- \mathbb{E} \widetilde{T}_{im})' H^{-1/2}a $.   We apply  Bernstein inequality condition on  $\{ T_{im}\}_{m=1}^{N_i}$. By elementary computations:
 \begin{align*}
{\rm var} (\mathcal{T}_1| \{T_{im}\}_{m=1}^{N_i}) & = \frac{1}{N_i^2} \sum_{m=1}^{N_i}\mathbb E \Big[\Big((\widetilde{T}_{im}- \mathbb{E} \widetilde{T}_{im})' H^{-1/2}a\Big)^2 \Big| a\Big]\notag\\
& =\frac {1}{N_i} \sum_{j=1}^p d_i^0(j) \Big( {a(j)}/{h_j^{1/2}} - (d_i^0)' H^{-1/2}a\Big)^2 \notag\\
& =  \frac{1}{N_i} \sum_{j=1}^p  \frac{d_i^0(j)}{h_j} a^2(j) - \frac{1}{N_i} \big[(d_i^0)' H^{-1/2} a\big]^2 \notag\\
& \leq  \Vert a\Vert^2/N_i, 
\end{align*}
where we used that fact $d_i^0(j) = e_j' Aw_i \leq e_j' A \mathbf{1}_K = h_j$.
Furthermore, with the individual bound $N^{-1} \max_t \{a(t)/\sqrt{h_t}\}$, we obtain from Bernstein inequality that with probability $1- o(n^{-(2C_0+4)})$:
\begin{align*}
|\mathcal{T}_1|\leq C\left( \sqrt{\frac{\log (n)} {N}} \, \Vert a\Vert + \frac 1N  \max_t \frac{|a(t)| }{\sqrt{h_t}} \log(n)  \right),
\end{align*}
by choosing appropriately large $C>0$. We then consider using Bernstein inequality to study $a(t)$ and get:
\[
|a(t)|\leq C \sqrt{\frac{\log (n)}{N}} + C\frac{\log (n)}{N\sqrt{h_{\min}}}
\]
simultaneously for all $1\leq t \leq p$, with probability $1- o(n^{-(2C_0+4)})$. As a result, under the condition $\min\{p, N\}\geq C_0\log(n)$ from Assumption~\ref{asm:para}, it holds that:
\begin{align} \label{bound-T1}
|\mathcal{T}_1| &\leq C\left( \sqrt{\frac{\log (n)} {N}} \, \Vert a\Vert + \frac 1N  \max_t \frac{|a(t)| }{\sqrt{h_t}} \log(n)  \right) \notag\\
& \leq C\left( \sqrt{\frac{p \log (n)} {N}} \,\Big[ \sqrt{\frac{\log (n)}{N}} + C\frac{\log (n)}{N\sqrt{h_{\min}}}\Big] + \frac pN\right) \notag\\
& \leq C\frac{p}{N}\, .
\end{align} 
We then proceed to  the second term in (\ref{form tI_1}), $\mathcal{T}_2 =N_i^{-2} \sum_{m=1}^{N_i}  (T_{im}- \mathbb{E} T_{im})' H^{-1} (\widetilde{T}_{im}- \mathbb{E} \widetilde{T}_{im}) $.  Using Bernstein inequality, similarly to the above derivations, we get:
\begin{align*}
{\rm var} (\mathcal{T}_2) &= N_i^{-4} \sum_{m=1}^{N_i} \mathbb E \Big( (T_{im}- \mathbb{E} T_{im})' H^{-1} (\widetilde{T}_{im}- \mathbb{E} \widetilde{T}_{im})\Big)^2\notag\\
& = N_i^{-4} \sum_{m=1}^{N_i} \mathbb E \bigg[\sum_{j=1}^p \frac{d_i^0(j)}{h_j^2} (\widetilde{T}_{im}(j)- \mathbb{E} \widetilde{T}_{im}(j))^2  - \Big((d_i^0)'H^{-1} (\widetilde{T}_{im}- \mathbb{E} \widetilde{T}_{im})\Big)^2 \bigg] \notag\\
& = N_i^{-3} \bigg[\sum_{j=1}^p \frac{(d_i^0(j))^2 (1 - d_i^0(j)) }{h_j^2}  -  \sum_{j=1}^p d_i^0(j) \Big( \frac{d_i^0(j)}{h_j} - (d_i^0)' H^{-1} d_i^0\Big)^2\bigg] \notag\\
& = N_i^{-3} \bigg[\sum_{j=1}^p \frac{(d_i^0(j))^2 (1 - 2d_i^0(j)) }{h_j^2}  + \Big( (d_i^0)' H^{-1} d_i^0\Big)^2\bigg] \notag\\
& < 2\frac{p}{N^3}\,. 
\end{align*}
The individual bound is given by  $N^{-2}/h_{\min}$. If follows from Bernstein inequality that:
 \begin{align}\label{22032302}
\mathcal{T}_2 \leq C\left(  \sqrt{\frac{p\log (n)}{N^3}} + \frac{\log (n)}{N^2h_{\min}}\right)
  \end{align}
 with probability $1- o(n^{-(2C_0+4)})$. Consequently, by pluging \eqref{bound-T1} and \eqref{22032302} into \eqref{form tI_1} and using Assumption~\ref{asm:para}, 
 \begin{equation} \label{eq:bdtI1}
  |\widetilde{\mathcal{I}}_1 | \lesssim  \frac pN
 \end{equation}
with probability $1- o(n^{-(2C_0+4)})$. By (\ref{eq:decoup}), we get:
\begin{align}\label{bd:I1}
  |{\mathcal{I}}_1 | \leq C\left( \sqrt{\frac{\log (n)} {N}} \, \Vert a\Vert   + \frac pN \right)
\end{align}
 with probability $1- o(n^{-(2C_0+4)})$.
 
Second, we prove a similar bound for $\mathcal{I}_2 $ with:
\[
 \mathcal{I}_2=  \frac{1}{N_i^2} \sum_{m=1}^{N_i} (T_{im}- \mathbb{E} T_{im})' H^{-1} (T_{im}- \mathbb{E} T_{im})  - \mathbb{E} (T_{im}- \mathbb{E} T_{im})' H^{-1} (T_{im}- \mathbb{E} T_{im}). 
\] 
We compute the variance by: 
  \begin{align*}
&\quad   {\rm var} (T_{im}- \mathbb{E} T_{im})' H^{-1} (T_{im}- \mathbb{E} T_{im})  \notag\\
   &= \mathbb{E} \Big( \sum_{t} h_t^{-1} (T_{im}(t) - d_i^0(t))^2\Big)^2 - \Big( \mathbb{E}  \sum_{t} h_t^{-1} (T_{im}(t) - d_i^0(t))^2\Big)^2 \notag\\
  &\leq  \sum_{t} h_t^{-2} d_i^0(t) \big[ (1- d_i^0(t))^4 + (1- d_i^0(t))d_i^0(t)^3\big] -  \sum_{t} h_t^{-2} d_i^0(t)^2 (1- d_i^0(t))^2 \notag\\
  &\leq \sum_{t} h_t^{-1} \lesssim ph_{\min}^{-1}.
  \end{align*}
This,   together with the crude bound: $$ |(T_{im}- \mathbb{E} T_{im})' H^{-1} (T_{im}- \mathbb{E} T_{im})  - \mathbb{E} (T_{im}- \mathbb{E} T_{im})' H^{-1} (T_{im}- \mathbb{E} T_{im})|\leq C h_{\min}^{-1}, $$
 gives  that with probability $1- o(n^{-(2C_0+4)})$, for some sufficiently large $C>0$:
 \begin{align}\label{bd:I2}
 |\mathcal{I}_2| \leq C\max\Big\{ \sqrt{\frac{p\log (n)}{N^3h_{\min}}}, \frac{\log (n)}{N^2h_{\min}}\Big\} \leq C \frac pN, 
 \end{align}
 under Assumption~\ref{asm:para}.
 Combing (\ref{bd:I1}) and (\ref{bd:I2}), yields that:
 \begin{align*}
 \tilde{z}_i' \tilde{z}_i  = z_i'H^{-1} z_i \leq  \mathbb{E} z_i' H^{-1} z_i +  |\mathcal{I}_1| +  |\mathcal{I}_2| \leq C\frac pN
 \end{align*}
 with probability $1- o(n^{-(2C_0+4)})$. Thus, we conclude the claim $ \mathbb{P} (E_i^c)=  o(n^{-(2C_0+4)})$ for all $1\leq i \leq n$. The proof of (a) is complete. 
 
%
%

Next, we show the proof of (b). 
Recall the second term on the RHS of (\ref{22040401}). Using the convexity of $\Vert \cdot \Vert$ and the trivial bound:
\[
\mathbb{E}| \tilde{z}_i' \tilde{z}_i  {\mathbf 1}_{E_i^c} |\leq \mathbb{P} ( {\mathcal{E}_i^c} ) \Vert  \tilde{z}_i' \tilde{z}_i \Vert_{\max} \leq h_{\min}^{-1} \mathbb{P} ( {\mathcal E_i^c} ),
\] 
we get:
\begin{align*}
\Big\Vert \sum_{i=1}^n \mathbb{E}  (\tilde{z}_i \tilde{z}_i' {\mathbf 1}_{\mathcal{E}_i^c} ) \Big\Vert \leq  \sum_{i=1}^n \mathbb{E}\big\Vert  \tilde{z}_i  \tilde{z}_i' {\mathbf 1}_{\mathcal{E}_i^c}\big\Vert =  \sum_{i=1}^n \mathbb{E}| \tilde{z}_i' \tilde{z}_i  {\mathbf 1}_{\mathcal{E}_i^c} | \leq  o(n^{-(2C_0+4)})  np =o(n^{-(C_0+3)})\,. 
\end{align*} 
Here, in the last step, we used the fact that $p \leq n^{C_0}$, which follows from the second condition in Assumption~\ref{asm:para}. 
This yields the estimate in  (b). 

Finally, we claim (c) by Matrix Bernstein inequality (i.e., Theorem~\ref{thm:matrixBernstein}). 
 Towards that, we need to derive the upper bounds of  $\Vert \widetilde{X}_i\Vert$ and $\Vert \mathbb{E} \widetilde{X}_i^2\Vert$. By definition of  $\widetilde{X}_i$, that is:
\[
\widetilde{X}_i := \frac 1n \big(\tilde{z}_i \tilde{z}_i' \mathbf{1}_{ \mathcal{E}_i} - \mathbb{E}  \tilde{z}_i \tilde{z}_i' \mathbf{1}_{\mathcal{E}_i} \big),
\]
we easily derive that:
\begin{align*}
\Vert \widetilde{X}_i\Vert \leq \frac 1n \Big( |\tilde{z}_i' \tilde{z}_i  {\mathbf 1}_{\mathcal{E}_i} | + \Vert \mathbb{E} (\tilde{z}_i \tilde{z}_i'  {\mathbf 1}_{\mathcal{E}_i}) \Vert \Big)
 \leq \frac1n \Big( |\tilde{z}_i' \tilde{z}_i  {\mathbf 1}_{\mathcal{E}_i} | + \Vert \mathbb{E} (\tilde{z}_i \tilde{z}_i'  {\mathbf 1}_{\mathcal{E}_i^c}) \Vert +\Vert \mathbb{E} (\tilde{z}_i \tilde{z}_i') \Vert   \Big)\leq\frac{Cp}{nN}
\end{align*}
for some large $C>0$, in which we used the estimate:
\begin{align*}
\Vert \mathbb{E} (\tilde{z}_i \tilde{z}_i') \Vert   = \Vert H^{-1/2} \mathbb{E} ( {z}_i {z}_i')H^{-1/2} \Vert  & \leq N_i^{-1} \Big\Vert H^{-1/2}\Big( {\rm diag} (d_i^0) - d_i^0 (d_i^0)'\Big)H^{-1/2} \Big\Vert \notag\\
& \leq  N_i^{-1} \Big\Vert H^{-1/2} {\rm diag} (d_i^0)H^{-1/2} \Big\Vert + N_i^{-1} \big|(d_i^0)'H^{-1}d_i^0\big| \notag\\
& \leq  \frac 2N \,. 
\end{align*}
By the above inequality, it also holds that:
\begin{align*}
\Vert \mathbb{E} (\tilde{z}_i \tilde{z}_i'  {\mathbf 1}_{\mathcal{E}_i}) \Vert  \leq  \Vert \mathbb{E} (\tilde{z}_i \tilde{z}_i'  {\mathbf 1}_{\mathcal{E}_i^c}) \Vert +\Vert \mathbb{E} (\tilde{z}_i \tilde{z}_i') \Vert
\leq \frac CN\, . 
\end{align*}
Moreover:
\begin{align*}
\Vert \mathbb{E} \widetilde{X}_i^2\Vert & = \big\Vert n^{-2} \mathbb E(\Vert \tilde z_i\Vert^2 \tilde z_i\tilde z_i' \mathbf{1}_{\mathcal{E}_i}) - n^{-2} (\mathbb{E} \tilde z_i\tilde z_i' \mathbf{1}_{\mathcal{E}_i})^2 \big\Vert \notag\\
&\leq \frac{p}{n^2N}  \Vert \mathbb E (\tilde z_i\tilde z_i' \mathbf{1}_{\mathcal{E}_i}) \Vert + \frac{1}{n^2} \Vert \mathbb E (\tilde z_i\tilde z_i' \mathbf{1}_{\mathcal{E}_i}) \Vert ^2 \notag\\
&  \leq \frac{Cp}{n^2N^2} \,. 
\end{align*}
Since  $\mathbb{E} \widetilde{X}_i=0 $, it follows from Theorem~\ref{thm:matrixBernstein} that:
\begin{align*}
\mathbb{P} \Big( \Big\Vert \sum_{i=1}^n \widetilde{X}_i  \Big\Vert \geq t \Big)\leq 2n \exp\Big( \frac{-t^2/2}{\sigma^2 + bt/3}\Big), 
\end{align*}
with $\sigma^2 = Cp/(nN^2)$, $b=Cp/(nN)$.  As a result: 
\begin{align*}
 \Big\Vert \sum_{i=1}^n \widetilde{X}_i  \Big\Vert\leq C \max\Big\{ \sqrt{\frac{p\log (n)}{nN^2}}, \frac{p\log (n)}{nN}\Big\}
\end{align*}
with probability $1- o(n^{-3})$, for some large $C>0$. We hence finish the proof of (c). The proof of \eqref{eq:tech2_E4} is complete now.

Lastly, we prove $\|E_1\| \leq C\sqrt{pn \log (n)}/\sqrt{N}$. By definition,  we rewrite:
\begin{align}
E_1  
& = (M^{-1/2} M_0^{1/2} ) M_0^{-1/2}DD' M_0^{-1/2} (M^{-1/2} M_0^{1/2} - I_p) \notag\\
& \quad  + (M^{-1/2} M_0^{1/2} - I_p) M_0^{-1/2}DD' M_0^{-1/2} \,. 
\end{align}
Decomposing $D$ by $ D_0 + Z$ gives rise to:
\begin{align}
\label{eq:MDDMU}
M_0^{-\frac 12}DD' M_0^{-\frac 12}& = M_0^{-\frac 12}\sum_{i=1}^n (1- N_i^{-1}) d_i^0 (d_i^0)' M_0^{-\frac 12}+ \frac nN I_p + M_0^{-\frac 12} D_0 Z' M_0^{-\frac 12}  + M_0^{-\frac 12} ZD_0' M_0^{-\frac 12} \notag\\
& \quad + M_0^{-\frac 12} (ZZ'- \mathbb{E} ZZ') M_0^{-\frac 12} \notag\\
& = G_0 + \frac nN I_p + E_2+ E_3+ E_4
\end{align}
Applying Lemma~\ref{lem:eigen}, together with (\ref{bd:E_2}) and (\ref{bd:new2}), we see that:
\begin{align*}
\Vert M_0^{-\frac 12}DD' M_0^{-\frac 12}\Vert \leq Cn
\end{align*}
Furthermore, it follows from  Lemma~\ref{lem:MM0} that:
\[
\|M^{-1/2} M_0^{1/2} - I_p\|  \leq C\sqrt{\frac{p\log (n)}{Nn}}, \qquad \text{and} \quad 
\|M^{-1/2} M_0^{1/2} \| = 1+o(1)\, . 
\]
Combining the estimates above, we conclude that:
\begin{align*}
\Vert E_1\Vert \leq C\sqrt{\frac{pn\log(n)}{N}}
\end{align*}
We therefore  finish the proof of Lemma~\ref{lem:tech1}. 

\subsection{Proof of Lemma~\ref{lem:tech2}}
We begin with the proof of (\ref{eq:entryE23}). Recall the definitions:
\[
E_2 = M_0^{-\frac 12} ZD_0' M_0^{-\frac 12}, \qquad E_3 = M_0^{-\frac 12} D_0Z' M_0^{-\frac 12}.  
\]
We bound:
\begin{align*}
\Vert e_j' E_2 \hat \Xi\Vert/n &\leq  \Vert e_j' E_2 \Vert/n \leq \frac{1}{n} \sum_{k=1}^K \Vert e_j'M_0^{-1/2} Z W_k\Vert  \cdot \Vert A_k'M_0^{-\frac 12}   \Vert \leq \frac{C}{n} \sum_{k=1}^K \Vert e_j'M_0^{-1/2} Z W_k\Vert 
\end{align*}
by the second inequality in (\ref{2024022302}). Similarly to how we derived (\ref{eq:M0ZW}), using Bernstein inequality, we have:
\begin{align*}
\Vert e_j'M_0^{-1/2} Z W_k\Vert  &\leq  C \frac{\sum_{i=1}^n z_i(j) W_k(i)}{\sqrt{h_j}} =  C{\sum_{i=1}^n \sum_{m=1}^{N_i} N_i^{-1}h_j^{-1/2}\big(T_{im} (j) - d_i^0(j)\big) W_k(i)} \notag\\
& \leq C\sqrt{\frac{\|W_k\|^2 \log (n)}{N}} + \frac{C\log (n)}{N\sqrt{h_{j}}} \notag\\
& \leq C\sqrt{\frac{n \log (n)}{N}} + \frac{C\log (n)}{N\sqrt{h_{j}}}
\end{align*}
with probability $1- o(n^{-C_0-3})$. Consequently:
\begin{align}\label{eq:e_jE_2Xi}
\Vert e_j' E_2 \hat \Xi\Vert/n \leq C\sqrt{ \frac{\log (n)}{Nn}} + C\frac{\log(n)}{nN\sqrt{h_j} }  \leq C\sqrt{ \frac{\log (n)}{Nn}} \leq C\sqrt{ \frac{h_j p\log (n)}{Nn}}
\end{align}
in view of $p\log (n)^2\leq Nn$ and $h_j \geq h_{\min} \geq c/p$ from Assumption~\ref{asm:para}.

Analogously, for $\Xi_3$, we have:
\begin{align}\label{eq:e_jE_3Xi}
& \Vert e_j' E_3 \hat \Xi \Vert /n \leq \frac{1}{n} \sum_{k=1}^K \Vert e_j' M_0^{-1/2} A_k \Vert \cdot \Vert W_k' Z' M_0^{-1/2} \hat \Xi\Vert \leq C\sqrt{\frac{h_j p\log (n)}{Nn}} \, .
\end{align}
where we used $\Vert W_k' Z' M_0^{-1/2} \hat \Xi\Vert\leq \Vert M_0^{-1/2} ZW_k\Vert\leq {\sqrt{pn\log (n)}}/{\sqrt{N}}$ from (\ref{eq:M0ZW}) and $\Vert e_j' M_0^{-1/2} A_k \Vert \leq C\sqrt{h_j}$.  Hence, we complete the proof of (\ref{eq:entryE23}).

In the sequel, we focus on the proof of (\ref{eq:tech2_entryE4}). Recall that $ E_4 = M_0^{-\frac12} (ZZ' - \mathbb{E} ZZ') M_0^{-\frac 12}$. We expect to show that:
\[
 \Vert e_j' E_4 \hat \Xi  \Vert /n  \leq  C\sqrt{\frac{h_j p\log (n)}{Nn}} \,  \Big(1 +  \Vert H_0^{-\frac 12} (\hat \Xi- \Xi O')\Vert_{2\to \infty} \Big) \,.
\]
Let us decompose $\Vert e_j' E_4 \hat \Xi  \Vert /n $ as follows:
\begin{align*}
n^{-1}\Vert e_j' E_4 \hat \Xi \Vert \leq n^{-1} \Vert e_j' E_4 \Xi\Vert  + n^{-1}\Vert e_j' E_4 (\hat \Xi - \Xi O')\Vert\,. 
\end{align*}
We bound $ n^{-1} \Vert e_j' E_4 \Xi\Vert$ first. 
For any fixed $1\leq k\leq K$, in light of the fact that $M_0(j,j) \asymp h_j$ for all $1\leq j\leq p$:
\begin{align*}
&| e_j' E_4 \xi_k| \asymp  |e_j' H^{-1/2} (ZZ' - \mathbb{E} ZZ') H^{-1/2} \xi_k| = \Big| \sum_{i=1}^n h_j^{-1/2} z_i(j)  z_i'H^{-1/2} \xi_k - h_j^{-1/2} \mathbb{E} z_i(j)  z_i'H^{-1/2}\xi_k\Big| \notag\\
&  =  \left| \sum_{i=1}^n \frac {1}{N_i^2} \sum_{m, m_1=1}^{N_i} \frac{T_{im}(j) - d_i^0(j)}{\sqrt{h_j}} \cdot (T_{im_1} - d_i^0)' H^{-\frac 12} \xi_k - \mathbb E \bigg[\frac{T_{im}(j) - d_i^0(j)}{\sqrt{h_j}} \cdot (T_{im_1} - d_i^0)' H^{-\frac 12} \xi_k \bigg]\right|  \notag\\
& \leq  |\mathcal{J}_1| + |\mathcal{J}_2|, 
\end{align*}
with: 
\begin{align*}
& \mathcal J_1 : =  \sum_{i=1}^n \frac {1}{N_i^2} \sum_{m}^{N_i} (T_{im} - d_i^0)' H^{-1/2}e_j \cdot (T_{im} - d_i^0)' H^{-1/2} \xi_k \notag\\
 & \qquad \qquad \qquad -\mathbb E (T_{im} - d_i^0)' H^{-1/2}e_j \cdot (T_{im} - d_i^0)' H^{-1/2} \xi_k, \notag\\
& \mathcal J_2 : = \sum_{i=1}^n \frac {1}{N_i^2} \sum_{m\neq m_1}^{N_i} (T_{im} - d_i^0)' H^{-1/2}e_j \cdot (T_{im_1} - d_i^0)' H^{-1/2} \xi_k\,. 
\end{align*}

For $ \mathcal J_1$,  it is easy to compute the order of its  variance as follows:
{\small
\begin{align*}
&\quad  {\rm var} ( \mathcal{J}_1)  \notag\\
& = \sum_{i=1}^n \sum_{m=1}^{N_i} N_i^{-4} {\rm var}  \Big( (T_{im} - d_i^0)' H^{-1/2}e_j \cdot (T_{im} - d_i^0)' H^{-1/2} \xi_k \Big)  \notag\\
& =  \sum_{i=1}^n \sum_{m=1}^{N_i} N_i^{-4} d_i^0(j) \cdot \frac{(1- d_i^0(j))^2}{h_j} \Big(  \frac{\xi_k(j)}{\sqrt{h_j}}  - \sum_{t} \frac{d_i^0(t)\xi_k(t)}{\sqrt{h_t}}\Big)^2  \notag\\
&+ \sum_{i=1}^n \sum_{m=1}^{N_i} N_i^{-4} \sum_{t\neq j} d_i^0(t)  \cdot \frac{( d_i^0(j))^2}{h_j} \Big(  \frac{\xi_k(t)}{\sqrt{h_t}}  - \sum_{s} \frac{d_i^0(s)\xi_k(s)}{\sqrt{h_s}}\Big)^2  \notag\\
& -  \sum_{i=1}^n \sum_{m=1}^{N_i}\frac{1}{N_i^4} \left(  \frac{d_i^0(j)}{\sqrt {h_j} }\Big(  \frac{\xi_k(j)}{\sqrt{h_j}}  - \sum_{t} \frac{d_i^0(t)\xi_k(t)}{\sqrt{h_t}}\Big) - \sum_{j=1}^p \frac{ (d_i^0(j))^2}{ \sqrt{h_j}} \Big(  \frac{\xi_k(j)}{\sqrt{h_j}}  - \sum_{t} \frac{d_i^0(t)\xi_k(t)}{\sqrt{h_t}}\Big) \right)^2
 \notag\\
& \leq C \frac{n}{N^3},
\end{align*}}where we used the facts that $\xi_k(t) \leq \sqrt{h_t}$, $d_i^0(j) \leq Ch_j$,  and $\sum_{t} d_i^0(t) = 1$. Furthermore, with the trivial bound of each summand in $\mathcal{J}_1$ given by $CN^{-2} h_j^{-1/2}$, it follows from the Bernstein inequality that:
\begin{align*}
|\mathcal{J}_1| \leq C\sqrt{\frac{n\log (n)}{N^3}} + C\frac{\log(n) }{N^2\sqrt{h_j}} \leq C\sqrt{\frac{n\log (n)}{N^3}}  
\end{align*}
with probability $1- o(n^{-3-C_0})$. Here, we used the conditions that $h_j \geq C/p$ and $p\log (n)^2 \leq Nn$.

We proceed to estimate $|\mathcal{J}_2|$. Employing Theorem~\ref{thm:decoupling2} with:
\[
h(T_{im}, T_{im_1}) = N_i^{-2} (T_{im} - d_i^0)' H^{-1/2}e_j \cdot (T_{im_1} - d_i^0)' H^{-1/2} \xi_k\,,
\]
 it suffices to examine the high probability bound of:
 \[
 \widetilde{\mathcal{J}}_2:= \sum_{i=1}^n \frac {1}{N_i^2} \sum_{m\neq m_1}^{N_i} (T_{im} - d_i^0)' H^{-1/2}e_j \cdot (\widetilde{T}_{im_1} - d_i^0)' H^{-1/2} \xi_k 
 \]
where $\{\widetilde{T}_{im_1}\}$ is an independent copy of $\{{T}_{im_1}\}$. Imitating the proof of (\ref{form tI_1}), we rewrite:
\[
\widetilde{\mathcal{J}}_2 = \sum_{i=1}^n \sum_{m=1}^{N_i} N_i^{-1} (T_{im} - d_i^0)' H^{-1/2}e_j \cdot b_{im}~~  \text{where}~~   b_{im}= \Big(\sum_{m_1\neq m} N_i^{-1} (\widetilde{T}_{im_1} - d_i^0)' H^{-1/2} \xi_k \Big)
\]
Notice that $b_{im}$ can be crudely bounded by $C$ in view of $\xi_k(t) \leq \sqrt{h_t}$.  Then, condition on $\{\widetilde T_{im_1}\}$, by Bernstein inequality, we can derive that:
\begin{align*}
|\widetilde{\mathcal{J}}_2| \leq C \Big( \sqrt{\frac{n \log (n)}{N}} + \frac{\log (n) }{N\sqrt{h_j}}\Big) \leq C  \sqrt{\frac{n \log (n)}{N}} 
\end{align*}
with probability $1- o(n^{-3-C_0})$. Consequently, we arrive at:
\begin{align*}
| e_j' E_4 \xi_k| \leq C  \sqrt{\frac{n \log (n)}{N}} \leq C  \sqrt{\frac{h_j pn \log (n)}{N}} 
\end{align*}
under the assumption that $h_j \geq C/p$. As $K$ is a fixed constant, we further conclude:
\begin{align}\label{est:ejE4U}
\Vert e_j' E_4 \Xi\Vert \leq C  \sqrt{\frac{h_j pn\log (n)}{N}}
\end{align}
with probability $1- o(n^{-3-C_0})$.

Next, we estimate $n^{-1}\Vert e_j' E_4 (\hat \Xi-\Xi O') \Vert$. By definition, we write:
\begin{align*}
\frac 1n\Vert e_j' E_4 (\hat \Xi-\Xi O')  \Vert = \frac 1n\Vert e_j' M_0^{-1/2} (ZZ'-\mathbb{E} ZZ' ) M_0^{-1/2}  (\hat \Xi-\Xi O') \Vert \,. 
\end{align*}
For each $1\leq t \leq p$:
\begin{align*}
&\quad \frac 1n | e_j' M_0^{-1/2} (ZZ'-\mathbb{E} ZZ' ) e_t | 
\notag\\
&\asymp \frac{1}{n\sqrt{h_j }}\sum_{i=1}^n z_i(j) z_i(t) - \mathbb{E} (z_i(j) z_j(t) ) \notag\\
& = \frac{1}{n\sqrt{h_j }}\sum_{i} \sum_{m, \tilde{m}}N_i^{-2} (T_{im}(j)- d_i^0(j)) (T_{i\tilde{m}}(t)- d_i^0(t))- \mathbb{E}(T_{im}(j)- d_i^0(j)) (T_{i\tilde{m}}(t)- d_i^0(t)) \notag\\
&=  \frac{1}{n \sqrt{h_j }}\sum_{i, m}N_i^{-2}(T_{im}(j)- d_i^0(j)) (T_{i{m}}(t)- d_i^0(t))- \mathbb{E}(T_{im}(j)- d_i^0(j)) (T_{i{m}}(t)- d_i^0(t)) \notag\\
& \quad  + \frac{1}{n \sqrt{h_j }}\sum_{i} N_i^{-2}\sum_{m\neq \tilde{m}} (T_{im}(j)- d_i^0(j)) (T_{i\tilde{m}}(t)- d_i^0(t))\notag\\
&:= (I)_t + (II)_t. 
\end{align*}

For $(I)_k$, using Bernstein inequality, it yields that with probability $1- o(n^{-3-2C_0})$:
\begin{align*}
\big| (I)_t \big| &\leq C \left\{
\begin{array}{lc}
 \max\Big\{ \sqrt{\frac{(h_j+h_t) h_t \log (n)}{nN^3}}, \, \frac{(h_j+ h_t)\log (n)}{nN^2\sqrt{h_j}}\Big\},  & t\neq j\\[0.2cm]
  \max\Big\{  \sqrt{\frac{ \log (n)}{nN^3}},\, \frac{\log (n)}{nN^2\sqrt{h_j}}  \Big\},  & t =j
\end{array}
\right. \notag\\
&\leq C \left\{
\begin{array}{lc}
  \sqrt{\frac{(h_j+h_t) h_t \log (n)}{nN^3}},   & t\neq j\\[0.2cm]
   \sqrt{\frac{ \log (n)}{nN^3}},  & t=j
\end{array}
\right.
\end{align*}
where the last step is due the the fact $p\log (n)^2 \leq Nn$ from Assumption~\ref{asm:para}.  
As a result:
\begin{align}\label{eq:l1_I}
\sum_{t=1}^p \big| (I)_t \big|\leq C\left(\sqrt {p}\sqrt{
\frac{\sum_{t\neq j} h_jh_t \log (n)}{nN^3}} + \sum_{t\neq j} h_t \sqrt{\frac{\log (n)}{nN^3}} + \sqrt{\frac{\log (n)}{nN^3}} \right)\leq C\sqrt{\frac{h_j p\log (n)}{nN^3}} 
\end{align}
Here, we used the Cauchy--Schwarz inequality to get:
\[
\sum_{t\neq j}  \sqrt{\frac{h_jh_t \log (n)}{nN^3}} \leq \sqrt{p-1}\cdot \sum_{t\neq j}  {\frac{h_jh_t \log (n)}{nN^3}} \leq \sqrt {p}\sqrt{
\frac{\sum_{t\neq j} h_jh_t \log (n)}{nN^3}} \,. 
\]

For $(II)_t$, since it is a U-statistics, we then apply the decoupling idea, i.e., Theorem~\ref{thm:decoupling2}, such that its high probability bound can be controlled by that of $(\widetilde{II})_t$, defined by:
\begin{align*}
(\widetilde{II})_t:= \frac{1}{n \sqrt{h_j }}\sum_{i} N_i^{-2}\sum_{m\neq \tilde{m}} ({T}_{im}(j)- d_i^0(j)) (\widetilde{T}_{i\tilde{m}}(t)- d_i^0(t)).
\end{align*}
where $\{\widetilde{T}_{i\tilde{m}}\}_{i, \tilde{m}}$ is the  i.i.d. copy of $\{{T}_{i {m}}\}_{i,m}$.  We further express: 
\begin{align*}
(\widetilde{II})_t &= \frac{1}{n \sqrt{h_j }}\sum_{i}N_i^{-2} \sum_{m} ({T}_{im}(j)- d_i^0(j)) {\bf \widetilde{T}}_{i, -m},
\end{align*}
where $ {\bf \widetilde{T}}_{i, -m} := \sum_{\tilde{m}\neq m} (\widetilde{T}_{i\tilde{m}}(t)- d_i^0(t))$. Condition on $\{\widetilde{T}_{i\tilde{m}}\}_{i, \tilde{m}}$, we use Bernstein inequality and get:
\begin{align*}
(\widetilde{II})_t &\leq C\max\Big\{ \sqrt{\frac{ \log (n) \cdot\sum_{i,m}  {\bf \widetilde{T}}_{i, -m}^2}{n^2N^4}}, \, \frac{\log (n)\cdot \max_{i,m}{|\bf \widetilde{T}}_{i, -m}| }{nN^2 \sqrt{h_j}}\Big\} \notag\\
& \leq C \sqrt{\frac{\log (n) \cdot  \max_{i,m}{|\bf \widetilde{T}}_{i, -m}|^2 }{nN^3}},
\end{align*}
in light of  $p\log (n)^2\leq Nn$. 
Furthermore, notice that:
\begin{align*}
\max_{i,m}{|\bf \widetilde{T}}_{i, -m}| \leq \sum_{\tilde{m}} \big| \widetilde{T}_{i\tilde{m}}(t)- d_i^0(t) \big|.
\end{align*}
It follows that:
\begin{align}\label{eq:l1_II}
\sum_{t=1}^ p \big| (\widetilde{II})_t \big| \leq C\sqrt{\frac{\log (n)}{nN}}\, \cdot \frac{1}{N}\sum_{t=1}^p\max_{i,m}{|\bf \widetilde{T}}_{i, -m}| & \leq C\sqrt{\frac{\log (n)}{nN}}\, \cdot \frac{1}{N}\sum_{t=1}^p \sum_{\tilde{m}} \big| \widetilde{T}_{i\tilde{m}}(t)- d_i^0(t) \big|\notag\\
&\leq C\sqrt{\frac{\log (n)}{nN}}\, , 
\end{align}
where the last step is due to the trivial bound that:
\[
\sum_{t=1}^p  \big| \widetilde{T}_{i\tilde{m}}(t)- d_i^0(t) \big| \leq 1+ \sum_{t=1}^p d_i^0(t) \leq C
\]
 for any $1\leq \tilde m \leq N$. 
Thus, combining (\ref{eq:l1_I}) and (\ref{eq:l1_II}), under  the condition $h_j\geq C/p$,  we obtain:
\begin{align} \label{2023101001}
 \frac 1n \|e_j' M_0^{-1/2} (ZZ'-\mathbb{E} ZZ' )\|_1 = \frac 1n\sum_{t=1}^p | e_j' M_0^{-1/2} (ZZ'-\mathbb{E} ZZ' ) e_t | \leq  C\sqrt{\frac{h_j p\log (n)}{nN}} 
\end{align}
with probability $1- o(n^{-3-C_0})$.

Moreover, employing the estimate $M_0(j,j) \asymp h_j$ for all $1\leq j \leq p$, it follows that:
\begin{align} \label{2022042902}
\frac 1n\Vert e_j' E_4 (\hat \Xi-\Xi O')  \Vert & = \frac 1n\Vert e_j' M_0^{-1/2} (ZZ'-\mathbb{E} ZZ' ) M_0^{-1/2}  (\hat \Xi-\Xi O') \Vert \notag\\
& \leq  \frac 1n \|e_j' M_0^{-1/2} (ZZ'-\mathbb{E} ZZ' )\|_1 \cdot \|M_0^{-1/2} H^{1/2} \| \cdot  \Vert H^{-1/2} (\hat \Xi- \Xi O')\Vert_{2\to \infty}\notag\\
&\leq C\sqrt{\frac{h_j p \log (n)}{nN}} \, \Vert H^{-1/2} (\hat \Xi- \Xi O')\Vert_{2\to \infty}\, 
\end{align}
with probability $1- o(n^{-3-C_0})$.

In the end, we combine (\ref{est:ejE4U}) and (\ref{2022042902}) and consider all $j$ simultaneously to conclude that:
\begin{align*}
n^{-1}\Vert e_j' E_4 \hat \Xi \Vert &\leq n^{-1} \Vert e_j' E_4 \Xi\Vert  + n^{-1}\Vert e_j' E_4 (\hat \Xi - \Xi O')\Vert   \\
&\leq  C\sqrt{\frac{h_j p \log (n)}{nN}} \,\Big(1+ \Vert H^{-1/2} (\hat \Xi- \Xi O')\Vert_{2\to \infty}\Big) \,
\end{align*}
with probability $1- o(n^{-3-C_0})$.  Combining all $1\leq j \leq p$, together with $p\leq n^{C_0}$, we complete the proof.

\subsection{Proof of Lemma~\ref{lem:tech3}}
We first prove (\ref{eq:tech2_add}) that:
\[
 \Vert e_j' E_4(M_0^{1/2}M^{-1/2}-I_p) \hat \Xi  \Vert/n \leq C\sqrt{h_j} \cdot \frac{p\log (n)}{nN} \Big( 1 +   \Vert H^{-\frac 12}  (\hat \Xi - \Xi O') \Vert_{2\to \infty}\Big) 
 \]
By the definition that $E_4 = M_0^{-1/2} (ZZ'-\mathbb{E} ZZ' ) M_0^{-1/2}$, we bound:
 \begin{align*}
\Vert e_j' E_4(M_0^{1/2}M^{-1/2}-I_p) \hat \Xi  \Vert/n  \leq \frac1n  \Vert e_j' M_0^{-1/2} (ZZ'-\mathbb{E} ZZ' ) \Vert_1 \cdot \Vert M_0^{-1/2} (M_0^{1/2}M^{-1/2}-I_p) \hat \Xi \Vert_{2\to \infty}\,. 
\end{align*}
From (\ref{2023101001}), it holds that $\Vert e_j' M_0^{-1/2} (ZZ'-\mathbb{E} ZZ' ) \Vert_1/n \leq C\sqrt{h_j p\log (n) }/\sqrt{nN} $ with probability $1- o(n^{-3-C_0})$. 
Next, we bound:
\begin{align*}
\Vert M_0^{-1/2} (M_0^{1/2}M^{-1/2}-I_p) \hat \Xi\Vert_{2\to \infty} &\leq \Vert H^{-1/2} (M_0^{1/2}M^{-1/2}-I_p) \Xi \Vert_{2\to \infty} 
\notag\\
&\quad + \Vert H^{-1/2} (M_0^{1/2}M^{-1/2}-I_p) (\hat \Xi-\Xi O') \Vert_{2\to \infty}
\end{align*}
The first term on the RHS can be bounded  simply by:
\begin{align*}
\Vert H^{-1/2} (M_0^{1/2}M^{-1/2}-I_p) \Xi \Vert_{2\to \infty} & \leq C\max_{i} |h_i^{-1/2} \sqrt{p\log (n)/nN} \cdot \sqrt {h_i}| \notag\\
& \leq C\sqrt{p\log (n)/nN} = o(1)
\end{align*}
The second term can be simplified to:
\begin{align*}
\Vert H^{-1/2} (M_0^{1/2}M^{-1/2}-I_p) (\hat \Xi-\Xi O') \Vert_{2\to \infty}&  
=\Vert  (M_0^{1/2}M^{-1/2}-I_p)  H^{-1/2}(\hat \Xi-\Xi O') \Vert_{2\to \infty} \notag\\
 &\leq C \sqrt{\frac{p\log (n)}{nN}}\cdot \Vert H^{-1/2}  (\hat \Xi-\Xi O') \Vert_{2\to \infty}.
\end{align*}
As a result:
\begin{align}\label{2022042901}
\Vert e_j' E_4(M_0^{1/2}M^{-1/2}-I_p) \hat \Xi  \Vert/n  & \leq C\sqrt{\frac{h_j p\log (n)}{nN}}\cdot\sqrt{\frac{p\log (n)}{nN}} \Big(1+    \Vert H_0^{-\frac 12}  (\Xi - \Xi_0O') \Vert_{2\to \infty}\Big)\notag\\
&  \leq C\sqrt{h_j} \cdot \frac{p\log (n)}{nN} \Big( 1 +   \Vert H^{-\frac 12}  (\hat \Xi - \Xi O') \Vert_{2\to \infty}\Big) 
\,. 
\end{align}
This proves \eqref{eq:tech2_add}. 

Subsequently, we prove (\ref{eq:tech3_add}) that:
\[
\Big\Vert e_j'  \big(M^{1/2}M_0^{-1/2} - I_p\big)\hat\Xi  \Big \Vert \leq C \sqrt{\frac{\log (n)}{Nn} } + o(\beta_n) \cdot \Vert e_j'(\hat \Xi - \Xi O')\Vert.
\]
We first bound:
\begin{align*}
\Big\Vert e_j'  \big(M^{1/2}M_0^{-1/2} - I_p\big)\hat\Xi  \Big \Vert& \leq  \Big\Vert e_j'  \big(M^{1/2}M_0^{-1/2} - I_p\big)\Xi  \Big \Vert + \Big\Vert e_j'  \big(M^{1/2}M_0^{-1/2} - I_p\big)(\hat\Xi - \Xi O') \Big \Vert.
\end{align*}
By Lemma~\ref{lem:MM0}, $|{M(j,j) - M_0(j,j)}|/{M_0(j,j)} \leq C\sqrt{\log(n)} /\sqrt{Nnh_j} $. It follows that:
\begin{align*}
 \Big\Vert e_j'  \big(M^{1/2}M_0^{-1/2} - I_p\big)\Xi  \Big \Vert &\leq \left|\sqrt{\frac{M(j,j)}{M_0(j,j)}} - 1\right| \cdot \|e_j'\Xi\| \notag\\
 & \leq C \frac{|{M(j,j) - M_0(j,j)}|}{M_0(j,j)} \cdot \|e_j'\Xi\| \notag\\
 & \leq C\sqrt{\frac{\log(n)}{Nn}}\,,
\end{align*}
and:
\begin{align*}
\Big\Vert e_j'  \big(M^{1/2}M_0^{-1/2} - I_p\big)(\hat\Xi - \Xi O') \Big \Vert  
&\leq  \left|\sqrt{\frac{M(j,j)}{M_0(j,j)}} - 1\right| \cdot \|e_j'(\hat \Xi - \Xi O')\| \notag\\
& \leq \sqrt{\frac{p\log(n)}{Nn}} \cdot  \|e_j'(\hat \Xi - \Xi O')\|\notag\\
& = o(\beta_n) \cdot \|e_j'(\hat \Xi - \Xi O')\| \, . 
\end{align*}
by the condition that $p \log(n) \ll Nn$. 
We therefore conclude (\ref{eq:tech3_add}), simultaneously for all $1\leq j \leq p$, with probability $1- o(n^{-3})$.

Lastly, we prove (\ref{eq:tech3_entryE1}). By the definition:
\[
E_1=  M^{-\frac 12} DD' M^{-\frac 12} - M_0^{-\frac 12} DD' M_0^{-\frac 12}, 
\]
and the decomposition:
\begin{align*}
M_0^{-\frac 12}DD' M_0^{-\frac 12}& =  G_0 + \frac nN I_p + E_2+ E_3+ E_4,  \text{ where } G_0 = M_0^{-1/2} \sum_{i=1}^n (1 - N_i^{-1} ) d_i^0 (d_i^0)' M_0^{-1/2}, 
\end{align*}
we bound:
\begin{align*}
&\quad \Vert e_j E_1 \hat \Xi \Vert/n   \notag\\
& \leq \Vert e_j' (I_p- M_0^{-1/2} M^{1/2} ) M^{-1/2} DD' M^{-1/2} \hat\Xi \Vert /n + \Vert  e_j' M_0^{-1/2} DD' M_0^{-1/2} (M_0^{1/2} M^{-1/2}-I_p) \hat\Xi   \Vert /n\notag\\
& \leq C\Vert e_j' (I_p- M_0^{-1/2} M^{1/2} )  \hat\Xi \Vert + C\Vert e_j' G_0 (M_0^{1/2}M^{-1/2}-I_p)\hat  \Xi  \Vert /n  \notag\\
& \qquad  + \Vert  e_j' (M_0^{1/2} M^{-1/2}-I_p) \hat\Xi   \Vert /N + \sum_{i=2}^4  \Vert e_j' E_i (M_0^{1/2}M^{-1/2}-I_p) \hat \Xi  \Vert/n,  
\end{align*}
where we used the fact that $M^{-1/2} DD' M^{-1/2} \hat\Xi = \widetilde \Lambda \hat\Xi$, where $\widetilde \Lambda =\hat  \Lambda + nN^{-1} I_p$, which leads to $\| \widetilde \Lambda\|\leq Cn$.

In the same manner to prove $\Vert e_j' E_2 \hat \Xi  \Vert/n$ and $\Vert e_j' E_3 \hat \Xi\Vert/n$, we can bound:
\begin{align} \label{2024022401}
\frac 1n\Vert e_j' E_s (M_0^{1/2}M^{-1/2}-I_p) \hat \Xi  \Vert \leq \frac 1n\Vert e_j' E_s \Vert \Vert M_0^{1/2}M^{-1/2}-I_p\Vert \leq  C\sqrt{\frac{h_jp  \log (n)}{Nn}}, \qquad \text { for } s=2,3\, . 
\end{align}
By Lemma~\ref{lem:MM0}, we derive:
\begin{align} \label{2024022402}
\Vert e_j' G_0 (M_0^{1/2}M^{-1/2}-I_p) \hat \Xi  \Vert /n&\leq C\sum_{t=1}^p  \frac{1}{\sqrt{h_jh_t}} |a_j' \Sigma_W a_t| \sqrt{\frac{ \log (n)}{h_tNn }} \Vert e_t ' \hat \Xi\Vert  \notag\\
& \leq C\sqrt{\frac{h_j p \log (n)}{Nn}}, 
\end{align}
where we crudely  bound $|a_j' \Sigma_W a_t| \leq  h_jh_t$, and use Cauchy--Schwarz inequality that $\sum_{t=1}^p \Vert e_t'  \hat \Xi\Vert \leq \sqrt{p} \sqrt{{\rm tr}(\hat \Xi \hat \Xi')}\leq K\sqrt p $. 
In addition:
\begin{align*}
 \Vert  e_j' (M_0^{1/2} M^{-1/2}-I_p) \hat\Xi   \Vert /N & \leq \left|\sqrt{M_0(j,j)}/\sqrt{M(j,j)}\right| \cdot \Vert e_j' (I_p- M_0^{-1/2} M^{1/2} )  \hat\Xi \Vert \notag\\
 & \leq C  \Vert e_j' (I_p- M_0^{-1/2} M^{1/2} )  \hat\Xi \Vert\,, 
\end{align*}
which results in:
\begin{align} \label{lem-A6-quote}
\Vert e_j E_1 \hat \Xi \Vert/n 
& \leq C\Vert e_j' (I_p- M_0^{-1/2} M^{1/2} )  \hat\Xi \Vert + C\Vert e_j' G_0 (M_0^{1/2}M^{-1/2}-I_p)\hat  \Xi  \Vert /n  \notag\\
& \qquad  + \sum_{i=2}^4  \Vert e_j' E_i (M_0^{1/2}M^{-1/2}-I_p) \hat \Xi  \Vert/n. 
\end{align}
Combining (\ref{2024022401}), (\ref{2024022402}), (\ref{eq:tech2_add}), and (\ref{eq:tech3_add}) into the above inequality, we  complete the proof of (\ref{eq:tech3_entryE1}).

\section{Proofs of the Rates for Topic Modeling}\label{appE}
The proofs in this section are quite similar to those  in \cite{ke2022using} by employing  the bounds in Theorem~\ref{thm:row_eigenv}. For readers' convenience, we provide brief sketches and refer to more details in the supplementary materials of \cite{ke2022using}. Notice that  $N_i \asymp \bar N \asymp N$ from Assumption~\ref{asm:para}. Therefore, throughout this section, we always assume $\bar{N}=N$ without loss of generality.

\subsection{Proof of Theorem~\ref{thm:row_R} }
Recall that:
\[
\hat R  = (\hat r_1, \hat r_2, \ldots, \hat r_p)' = [{\rm diag}(\hat \xi_1)]^{-1} (\hat \xi_2, \ldots, \xi_K).
\]
Since the first eigenvector of $G_0$ is with multiplicity one, which can been seen in Lemma~\ref{lem:eigen}, and the fact that $\|G - G_0\|\ll n$, it is not hard to obtain that $O' = {\rm diag} (\omega, \Omega')$ where $\omega \in \{1, -1\}$ and $ \Omega'$ is an orthogonal matrix in $\mathbb R^{K-1, K-1}$. Let us write $\hat \Xi_1 : = (\hat \xi_2, \ldots, \hat \xi_K)$ and similarly for $\Xi_1$. Without loss of generality, we assume $\omega =1$. Therefore:
\begin{align} \label{2024022501}
\big| \xi_1(j) - \hat \xi_1(j)  \big|  \leq C\sqrt{\frac{h_j p\log(n)}{Nn\beta_n^2}}, \qquad \big\| e_j'(\hat \Xi_1 - \Xi_1) \Omega'  \big\|  \leq C\sqrt{\frac{h_j p\log(n)}{Nn\beta_n^2}} \,. 
\end{align}
We rewrite:
\begin{align*}
\hat r_j' - r_j' \Omega' =\hat \Xi_1(j) \cdot  \frac{\xi_1(j) - \hat \xi_1(j) }{\hat \xi_1(j)\xi_1(j)}  -  \frac{e_j' (\hat \Xi_1 - \Xi_1 \Omega' ) }{ \xi_1(j)}. 
\end{align*}
Using Lemma~\ref{lem:eigen} together with (\ref{2024022501}), we conclude the proof. 

\subsection{Proof of Theorem~\ref{thm:A_k}}
In this section, we provide a simplified proof by neglecting the details about some quantities in the oracle case. We refer readers to the proof of Theorem 3.3 of  \cite{ke2022using}  for more rigorous arguments. 

Recall the Topic-SCORE algorithm. Let $\widehat V = (\hat v_1, \hat v_2, \ldots, \hat v_K )$ and denote its population counterpart by $V$. We write:
\begin{align*}
\hat Q = \left(
\begin{array}{ccc}
1& \ldots  &1 \\
\hat v_1 & \ldots & \hat v_K
\end{array}
\right), 
\qquad 
Q = \left(
\begin{array}{ccc}
1& \ldots  &1 \\
 v_1 & \ldots & v_K
\end{array}
\right)
\end{align*}
Similarly to  \cite{ke2022using}, by properly choosing the vertex hunting algorithm and the anchor words condition, it can be seen that:
\[
\| \hat V   -  V \| \leq C\sqrt{\frac{p\log(n)}{Nn\beta_n^2}}
\]
where we omit the permutation for simplicity here and throughout this proof. 
As a result:
\begin{align*}
\| \hat \pi_j^*  - \pi_j^*\| & =\left\| \hat Q^{-1} \left(\begin{array}{c}1\\ \hat r_j\end{array} \right) -
 Q^{-1} \left(\begin{array}{c}1\\  \Omega r_j\end{array} \right)  \right\| \notag\\
 & \leq C \|Q^{-1}\|^2 \cdot \| \hat V   -  V \| \cdot \|r_j\|  + \|Q^{-1} \| \|\hat r_j - \Omega r_j\|\notag\\
 & \leq C\sqrt{\frac{p\log(n)}{Nn\beta_n^2}} = o(1)
\end{align*}
where we used the fact that  $\|Q^{-1}\|\leq C$  whose details can be found in the proof of Lemma~G.1 in supplementary material of \cite{ke2022using}. Considering the truncation at 0, it is not hard to see that:
\[
\|\tilde \pi_j^* - \pi_j^*\| \leq C \| \hat \pi_j^*  - \pi_j^*\|\leq C\sqrt{\frac{p\log(n)}{Nn\beta_n^2}} =o(1);
\]
and furthermore:
\begin{align}\label{2024022601}
\|\hat \pi_j - \pi_j \|_1 & \leq \frac{\|\tilde \pi_j^* - \pi_j^*\|_1}{\|\tilde \pi_j^*\|_1} + \frac{\|\pi_j^*\|_1 \big|\|\tilde \pi_j^*\|_1 -\| \pi_j^*\|_1\big|}{\|\tilde \pi_j^*\|_1\| \pi_j^*\|_1}\notag\\
& \leq C \|\tilde \pi_j^* - \pi_j^*\|_1 \leq C\sqrt{\frac{p\log(n)}{Nn\beta_n^2}}. 
\end{align}
by noticing that $\pi_j = \pi_j^*$ in the oracle case. 

Recall that $\tilde A = M^{1/2} {\rm diag}(\hat \xi_1)\hat \Pi =:(\tilde a_1, \ldots, \tilde a_p)'$. 
Let $A^* = M_0^{1/2}{\rm diag} (\xi_1) \Pi = (a^*_1, \ldots, a^*_p)'$. Note that  $A = A^* [{\rm diag}({\bf 1}_p A^*)]^{-1}$. We can derive:
\begin{align}\label{2024022610}
\|\tilde a_j - a^*_j \|_1&  \leq  \Big\|  \sqrt{M(j,j)}  \, \hat \xi_1(j)\hat \pi_j  -  \sqrt{M_0(j,j)}  \, \xi_1(j) \pi_j \Big\|_1 \notag\\
& \leq C \|\sqrt{M(j,j)}  - \sqrt{M_0(j,j)} \| \cdot \| \xi_1(j)\| \cdot \|\pi_j\|_1 + C\sqrt{M_0(j,j)}\, \cdot  \| \hat \xi_1(j) - \xi_1(j)\|\cdot \|\pi_j\|_1  \notag\\
&\qquad + C\sqrt{M_0(j,j)}\, \cdot  \| \xi_1(j)\|\cdot \|\hat \pi_j - \pi_j\|_1  \notag\\
& \leq C h_j \sqrt{\frac{p \log (n) }{Nn\beta_n^2}},
\end{align}
where we used (\ref{2024022601}), (\ref{2024022501}) and also Lemma~\ref{lem:MM0}.  Write $\tilde A = (\tilde A_1, \ldots, \tilde A_K)$ and $A^* = (A_1^*, \ldots, A_K^*)$.   We crudely bound: 
\begin{align}\label{2024022611}
\Big| \|\tilde A_k\|_1  - \|A^*_k\|_1\Big|\leq  \sum_{j=1}^p \|\tilde a_j - a^*_j \|_1 \leq
C \sqrt{\frac{p \log (n) }{Nn\beta_n^2}} =o(1) 
\end{align}
simultaneously for all $1\leq k \leq K$,
since $\sum_j h_j = K$. By the study of oracle case in \cite{ke2022using}, it can be deduced that $\|A^*_k\|_1\asymp 1$ (see more details in the supplementary materials of  \cite{ke2022using}). 
It then follows that:
\begin{align*}
\|\hat a_j - a_j \|_1 &=\Big\|{\rm diag} (1/\|\tilde A_1\|_1, \ldots, 1/\|\tilde A_K\|_1 )  \tilde a_j  - {\rm diag} (1/\| A^*_1\|_1, \ldots, 1/\|A^*_K\|_1 )  a^*_j \Big\|_1 \notag\\
& = \sum_{k=1}^K \bigg|  \frac{\tilde a_j(k) }{\|\tilde A_k\|_1} -  \frac{a^*_j(k)}{\|A^*_k\|_1} \bigg|\notag\\
& \leq  \sum_{k=1}^K \bigg|  \frac{\tilde a_j(k) - a^*_j(k) }{\|A^*_k\|_1} \bigg| + |a^*_j(k)| \frac{\big|\|\hat A_k\|_1 - \|A^*_k\|_1 \big|}{\|A^*_k\|_1\|\tilde A_k\|_1} \notag\\
& \leq C \sum_{k=1}^K \|\tilde a_j - a^*_j \|_1   + \|a_j^*\|_1 \max_{k}\Big| \|\tilde A_k\|_1  - \|A^*_k\|_1\Big| \notag\\
&  \leq
Ch_j  \sqrt{\frac{p \log (n) }{Nn\beta_n^2}}  = C\|a_j\|_1 \sqrt{\frac{p \log (n) }{Nn\beta_n^2}} \,. 
\end{align*}
Here, we used (\ref{2024022610}), (\ref{2024022611}) and  the following  estimate:
\[
\|a_j^*\|_1 = \sqrt{M_0(j,j)}\,   |\xi_1(j)| \|\pi^*_j\|\asymp h_j 
\]
Combining all $j$ together, we immediately have the result for $\mathcal L (\hat A, A)$.

\subsection{Proof of Theorem~\ref{thm:w_i}}
The optimization in (\ref{estimateW}) has a explicit solution given by:
\[
\hat w^*_i = \big(\hat A' M^{-1} \hat A\big)^{-1}  \hat A' M^{-1} d_i\,. 
\]
Notice that $ (A'M_0^{-1} A)^{-1} A'M_0^{-1}d^0_i =(A'M_0^{-1} A)^{-1} A'M_0^{-1}Aw_i = w_i $. 
Consequently:
\begin{align}\label{eq:hwi*-wi}
\|\hat w_i^* - w_i \|_1&  = \big\|  \big(\hat A' M^{-1} \hat A\big)^{-1}  \hat A' M^{-1} d_i  -  (A'M_0^{-1} A)^{-1} A'M_0^{-1}d^0_i \big\|_1 \notag\\
& \leq  \big\| (A'M_0^{-1} A)^{-1} \big( \hat A' M^{-1} \hat A -  A'M_0^{-1} A\big) \big(\hat A' M^{-1} \hat A\big)^{-1}  \hat A' M^{-1} d_i  \big\|_1\notag\\
& \quad + \big\| (A'M_0^{-1} A)^{-1}\big( \hat A' M^{-1} d_i -  A'M_0^{-1}d^0_i \big)  \big\|_1\notag\\
& \leq C\beta_n^{-1}  \| \big( \hat A' M^{-1} \hat A -  A'M_0^{-1} A \big) \| (\|\hat w^*_i - w_i \|_1 + \|w_i \|_1 )  \notag\\
& \quad + C \beta_n^{-1} \big\| \hat A' M^{-1} d_i -  A'M_0^{-1}d^0_i   \big\|, 
\end{align}
since $\big\|(A'M_0^{-1} A)^{-1}\big\|\asymp \big\|(A'H^{-1} A)^{-1}\big\| \asymp 1$.
What remains is to analyze:
\[
T_1 : =  \| \big( \hat A' M^{-1} \hat A -  A'M_0^{-1} A \big)  \|, \quad \text{and} \quad T_2: = \big\| \hat A' M^{-1} d_i -  A'M_0^{-1}d^0_i   \big\|.
\]

For $T_1$, we bound:
\begin{align*}
T_1 &\leq \| ( \hat A - A)' M^{-1} \hat A \| + \| A' ( M^{-1}- M_0^{-1}) \hat A \| \notag\\
& \quad+  \| A'M_0^{-1} (\hat A-A)\| \,. 
\end{align*} 
Using the estimates:
\begin{align*}
\| \hat a_j - a_j  \|_1   \leq C h_j \sqrt{\frac{p\log(n)}{Nn\beta_n^2}}, \qquad  
\big|M(j,j)^{-1}- M_0(j,j) ^{-1}\big| \leq \frac{\sqrt{ \log(n) }}{h_j\sqrt{Nnh_j}},
\end{align*}
it follows that:
\begin{align*}
\|A'(M^{-1}- M_0^{-1} ) (\hat A - A) \| &\leq \sum_{k,k_1=1}^K  \big|A_k(M^{-1}- M_0^{-1} )(\hat A_{k_1} - A_{k_1})  \big| \notag\\
& \ll  \sum_{k=1}^K  \|\hat A_{k} - A_{k}  \|_1 = \sum_{j=1}^p \| \hat a_j - a_j  \|_1 \notag\\
& \ll \sqrt{\frac{p\log(n)}{Nn\beta_n^2}}\,, 
\end{align*}
and similarly:
\begin{align*}
&\| (\hat A - A)' M_0^{-1}  (\hat A - A) \| 
  \ll  \sum_{k=1}^K  \|\hat A_{k} - A_{k}  \|_1 \ll \sqrt{\frac{p\log(n)}{Nn\beta_n^2}}, \notag\\
  & \| (\hat A - A)'(M^{-1}- M_0^{-1} ) (\hat A - A) \| \ll  \sum_{k=1}^K  \|\hat A_{k} - A_{k}  \|_1 \ll \sqrt{\frac{p\log(n)}{Nn\beta_n^2}}\,. 
\end{align*}
As a result:
\begin{align}\label{eq:24T1}
T_1 & \leq  C\|( \hat A - A)' M_0^{-1} A \| + C \| A' ( M^{-1}- M_0^{-1}) A \|\notag\\
& \leq C\sum_{j=1}^p \|\hat a_j - a_j\|_1 + C \sqrt{\frac{p\log(n)}{Nn}} \cdot \sum_{j=1}^p \|a_j\|_1 \notag\\
& \leq C  \sqrt{\frac{p\log(n)}{Nn\beta_n^2}}\,. 
\end{align}

Next, for $T_2$, we bound:
\begin{align*}
T_2 &\leq \| (\hat A - A)' M^{-1} d_i \| +  \| A'(M^{-1}- M_0^{-1} )  d_i \| + \| A'M_0^{-1} (d_i - d_i^0) \| \notag\\
& \leq\max_{j} \left( \frac{\|\hat a_j - a_j\|_1 }{h_j}  + \|a_j\|_1\frac{\sqrt{\log(n)}}{h_j\sqrt{Nnh_j}} \right)\cdot \|d_i\|_1 + \max_{1\leq k\leq K} \big|A_k'M_0^{-1}(d_i - d_i^0) \big|\notag\\
& \leq  C  \sqrt{\frac{p\log(n)}{Nn\beta_n^2}} + \max_{1\leq k\leq K} \big|A_k'M_0^{-1}(d_i - d_i^0) \big|\,. 
\end{align*}
where for $ (\hat A - A)' M^{-1} d_i $, given the low-dimension $K$,  we  crudely bound:
\[
\| (\hat A - A)' M^{-1} d_i \| \leq C\max_k \big| (\hat A_k - A_k)' M^{-1} d_i \big|\leq C\max_{k,j}\big| h_j^{-1}\big(\hat a_j(k) - a_j(k)\big)   \big| \|d_i\|_1
\]
and $\big|\hat a_j(k) - a_j(k)\big|\leq \|\hat a_j - a_j\|_1$. We bound $ \| A'(M^{-1}- M_0^{-1} )  d_i \|$ in the same manner. 
To proceed, we analyze $\big|A_k'M_0^{-1}(d_i - d_i^0) \big|$ for a fixed $k$. We rewrite it as:
\begin{align*}
A_k'M_0^{-1}(d_i - d_i^0) =\frac{1}{N_i} \sum_{m=1}^{N_i}A_k'M_0^{-1}( T_{im} - \mathbb T_{im}). 
\end{align*}
The RHS is an independent sum where Bernstein inequality can be applied. By elementary computations, the variance is:
\begin{align*}
N_i^{-1} {\rm var} \big( A_k'M_0^{-1}( T_{im} - \mathbb T_{im}) \big) & = N_i^{-1} \mathbb E\big( A_k'M_0^{-1}( T_{im} - \mathbb T_{im} )\big)^2\notag\\
& = N_i^{-1} A_k'M_0^{-1}{\rm diag } (d_i^0) M_0^{-1} A_k - N_i^{-1}\big( A_k'M_0^{-1}d_i^0\big)^2 \notag\\
& \leq N^{-1} 
\end{align*}
and the individual bound is crudely $N^{-1}$. It follows from Bernstein inequality that with probability $1- o(n^{-4})$:
\[
\|A_k'M_0^{-1}(d_i - d_i^0) \|\leq C\Big( \sqrt{\frac{\log (n)}{N}} + \frac{\log (n)}{N}\Big) \leq C \sqrt{\frac{\log (n)}{N}} 
\]
in light of $N \gg \log(n)$. This gives rise to:
\begin{align*}
T_2  \leq  C  \sqrt{\frac{p\log(n)}{Nn\beta_n^2}} + C \sqrt{\frac{\log (n)}{N}} 
\end{align*}

We substitute the above equation, together with (\ref{eq:24T1}), into  (\ref{eq:hwi*-wi}) and conclude that: 
\[
\|\hat w_i^* - w_i \|_1\leq  C  \sqrt{\frac{p\log(n)}{Nn\beta_n^4}} + C \sqrt{\frac{\log (n)}{N\beta_n^2}} \,. 
\]
Recall that the actual estimator $\hat w_i$ is defined by:
\[
\hat w_i  = \max\{ \hat w_i^*, 0\}/ \| \max\{ \hat w_i^*, 0\}\|_1, 
\]
where the maximum is taken entry-wisely. We write $\tilde w_i : =  \max\{ \hat w_i^*, 0\}$ for short. Since $w_i$ is always non-negative, it is not hard to see that:
\[
\|\tilde w_i  - w_i \|_1\leq C\|\hat w_i^* - w_i \|_1\leq  C  \sqrt{\frac{p\log(n)}{Nn\beta_n^4}} + C \sqrt{\frac{\log (n)}{N\beta_n^2}} = o(1)\,.
\]
As a result, $\|\tilde w_i \|_1 = 1+ o(1)$. Moreover:
\begin{align*}
\|\hat w_i - w_i \|_1 &\leq \frac{\|\tilde w_i - w_i\|_1}{\|\tilde w_i \|_1} + \|w_i\|_1 \Big|\frac{1}{\|\tilde w_i \|_1} - \frac{1}{\| w_i \|_1} \Big| \notag\\
& \leq C\|\tilde w_i - w_i\|_1 \leq  C  \sqrt{\frac{p\log(n)}{Nn\beta_n^4}} + C \sqrt{\frac{\log (n)}{N\beta_n^2}}
\end{align*}
with probability $1- o(n^{-4})$. Combining all $i$, we thus conclude the proof.

\bibliography{topic}

\end{document}